\documentclass[]{article}

\usepackage{amsmath}
\usepackage{amsfonts}
\usepackage{cite}
\usepackage{tikz}
\newtheorem{mydef}{Definition}
\newtheorem{mylem}{Lemma}
\newtheorem{mythm}{Theorem}

\usepackage{geometry}
\title{Kruskal-Katona type problem}
\author{Matthew Fitch \footnote{Zeeman Building,
University of Warwick
Coventry CV4 7AL, United Kingdom. Research supported by ERC Grant No. 306493. Email address: M.H.D.Fitch@warwick.ac.uk}}

\begin{document}

\maketitle

\begin{abstract}
The Kruskal Katona theorem was proved in the 1960s ~\cite{kruskal, katona}. In the theorem, we are given an integer $r$ and families of sets $\mathcal{A}\subset \mathbb{N}^{(r)}$ and $\mathcal{B}\subset\mathbb{N}^{(r-1)}$ such that for every $A\in\mathcal{A}$, every subset of $A$ of size $r-1$ is in $\mathcal{B}$. We are interested in finding the mimimum size of $b=|\mathcal{B}|$ given fixed values of $r$ and $a=|\mathcal{A}|$. The Kruskal Katona theorem states that this mimimum occurs when both $\mathcal{A}$ and $\mathcal{B}$ are initial segments of the colexicographic ordering. The Kruskal Katona theorem is very useful and has had many applications and generalisations.

In this paper, we are interested in one particular generalisation, where instead of every subset of $A$ of size $r-1$ being in $\mathcal{B}$, we will instead ask that only $k$ of them are, where $k$ is some integer smaller than $r$. Note that setting $k=r$ is exactly the Kruskal Katona theorem. We will first find exact results for the cases where $0\leq k \leq 3$. For $k\geq 4$ we will not solve the question completely, however, we will find the exact result for infinitely many $a$. We will also provide a formula that is within some additive constant of the correct result for all $a$.
	
\end{abstract}

\section{Introduction}

The Kruskal Katona therom was proved in the 1960s by Kruskal and Katona \cite{kruskal} ~\cite{katona}. In the theorem, we have families $\mathcal{A}\subset \mathbb{N}^{(r)}$ and $\mathcal{B}\subset\mathbb{N}^{(r-1)}$ such that for every $A\in\mathcal{A}$, there exist distinct $B_1, B_2,.. B_r\in\mathcal{B}$ with $B_i\subset A$. Given $a=|\mathcal{A}|$, we want to minimise $b=|\mathcal{B}|$, or equivalently, given $b$, we want to maximise $a$. Throughout this paper, we will find it simpler to use the formulation where $b$ is known and $a$ is being maximised. The Kruskal Katona theorem states that an optimal solution is when $\mathcal{A}$ is an initial segment of the colexicographic ordering on sets of size $r$, and $\mathcal{B}$ is the corresponding initial segment of the colexicographic ordering on sets of size $r-1$. \\

\begin{mydef}[Colexicographic ordering or colex]
The colexicographic ordering of sets of size $r$ is a total ordering where $A < B$ if the largest element of $B\backslash A$ is larger than the largest element of $A\backslash B$.

\end{mydef}	

\textbf{Example:}
When $r=4$, the first sets in the colex ordering are:\\ $\{1,2,3,4\}, \{1,2,3,5\}, \{1,2,4,5\}, \{1,3,4,5\}, \{2,3,4,5\}, \{1,2,3,6\}, \{1,2,4,6\}, \{1,3,4,6\}, ...$ . \\
\\
More specifically, we can calculate what $a$ or $b$ is based on the other. Pick some intergers $c_{r-1} > c_{r-2} > c_{r-3} > ... > c_1$ in such a way as to make the following decomposition of $b$ into binomials hold: 

\[b=\binom{c_{r-1}}{r-1} + \binom{c_{r-2}}{r-2} + \binom{c_{r-3}}{r-3} +... + \binom{c_{1}}{1}.\]

Note that every integer $b$ does have a valid set of $c_i$s that satisfy this formula, because $\binom{x}{s}+\binom{x-1}{s-1}+\binom{x-2}{s-2}+...+\binom{x-s+1}{1}=\binom{x+1}{s}-1$. So we can construct this set of $c_i$s it by first picking $c_{r-1}$ such that satisfies $\binom{c_{r-1}}{r-1} \leq b < \binom{c_{r-1}+1}{r-1}$. We subtract it from $b$ and then repeat the operation to find $c_{r-2}$. The remainder is between 0 and $\binom{c_{r-1}}{r-2}-1$ so when we pick $c_{r-2}$, it will be strictly less than $c_{r-1}$. Repeat in the same way to find $c_{r-3}$, $c_{r-4}$, ... $c_1$. Therefore such a decomposition exists for all $b$. \\
As for uniqueness, we note that if we pick a different decomposition, there is some largest $i$ where we picked different values for $c_i$. Say we use $c_{i}+t$ instead of $c_i$. Then  $\binom{c_{r-1}}{r-1} + ... + \binom{c_i+t}{i}>b$ by definition of how we picked $c_i$ to be the maximum number that worked. If on the other hand, we use $c_i-t$ instead of $c_i$, then the maximum number we can get is $\binom{c_{r-1}}{r-1} + ... + \binom{c_i-t+1}{i}-1 < b$ so this also doesn't work. Therefore the choice of $c_i$s is unique.\\
\\
Once we have written $b$ in that form, then the maximum $a$ is given by:

\[\binom{c_{r-1}}{r} + \binom{c_{r-2}}{r-1} + \binom{c_{r-3}}{r-2} +... + \binom{c_{1}}{2}.\]

Conversely, given $a=\binom{c_{r-1}}{r} + \binom{c_{r-2}}{r-1} + \binom{c_{r-3}}{r-2} +... + \binom{c_{1}}{2}+\binom{c_0}{1}$, the minimum $b$ is $b=\binom{c_{r-1}}{r-1} + \binom{c_{r-2}}{r-2} + \binom{c_{r-3}}{r-3} +... + \binom{c_{1}}{1}+[1 \text{ if $c_0>0$}]$.\\
\\
\\
In 2015, Bollob\`{a}s and Eccles asked if the Kruskal Katona theorem could be generalised to where instead of every subset of $A$ being in $\mathcal{B}$, we instead required only $k$ out fo the $r$ possible. ~\cite{bollobas}. In this case, we call the maximum value for $a$ given $r,k$ and $b$ to be $f(r,k,b)$. They considered one configuration in particular: let $\mathcal{A}$ be of the form $\{S\cup X\}$ where $X$ runs over an initial segment of the colex on sets of size $k$ and $S$ is just some set of size $r-k$ (that doesn't intersect any of the $X$s). Meanwhile let $\mathcal{B}$ be defined in the same way as $\{S\cup Y\}$ where $Y$ runs over the corresponding initial segment of the colex on sets of size $k-1$. These collections of sets have the property that for every $A$ in $\mathcal{A}$, there exist $B_1,B_2,...,B_k$ in $\mathcal{B}$ with $B_i\subset A$.\\
\\
This example gives a very similar formula to the Kruskal Katona theorem. It shows that for any $c_{k-1}>c_{k-2}>...>c_1$,

\[f\left(r,k,\binom{c_{k-1}}{k-1} + \binom{c_{r-2}}{k-2} + \binom{c_{k-3}}{k-3} +... + \binom{c_{1}}{1}\right) \geq \binom{c_{r-1}}{k} + \binom{c_{k-2}}{k-1} + \binom{c_{k-3}}{k-2} +... + \binom{c_{1}}{2}.\]

Bollob\`{a}s and Eccles conjectured that this configuration is actually the optimal one when $a$ and $b$ are large. They did also note that this conjecture cannot be extended to small values of $a$ and $b$, because they found an example that shows that $f(5,4,13)=6$. If you tried to use the conjecture, it would tell you the answer is 7, which is incorrect. So the example is not optimal for small values of $a$ and $b$; however, they still think this example is optimal when $a$ and $b$ are large enough.\\
\\
In this paper, we will first start by doing the easy cases of $k=0$, $k=1$, $k=2$ and $k=3$, which are all done using similar methods, although it gets more complicated as $k$ increases.

\begin{mythm} (the cases where $k\leq 3$)\\
$\bullet$ For $0=k\leq r$, the optimal value is $b=0$ regardless of what $a$ is.\\
$\bullet$ For $1=k\leq r$, the optimal value is $b=1$ if $a\geq 1$, otherwise it is $b=0$ if $a=0$.\\
$\bullet$ For $2=k\leq r$,  $f(r,2,b)=\binom{b}{2}$.\\
$\bullet$ For $3=k\leq r$, $f\left(r,3,\binom{c_2}{2}+\binom{c_1}{1}\right)=\binom{c_2}{3}+\binom{c_1}{2}$ whenever $c_2>c_1$ and $c_2\geq 29$.
\end{mythm}

Note that this is still missing the cases where $c_2<29$; however there are only finitely many of these so they can in theory be solved by simply checking every single case. After this, we will move on to the case where $k\geq 4$, which continues to use the same method, although due to some new complexities, we can no longer find an exact result for all $b$ large enough. However, the method still gives exact results for an infinite number of values for $b$:

\begin{mythm}
Given $0\leq k\leq r$, there is some constant $\mu$ depending only on $k$ such that if $c_{k-1}>c_{k-2}>...>c_1>\mu$, then:

\[f\left(r,k,\binom{c_{k-1}}{k-1} + ... + \binom{c_{1}}{1}\right)=\binom{c_{k-1}}{k} + ... + \binom{c_{1}}{2}.\]

\end{mythm}

Our method also allows us to get to within some additive contant of the answer for all $b$:

\begin{mythm}
	Given $0\leq k\leq r$, there is a constant $\tau$ depending only on $k$ such that if $b=\binom{c_{k-1}}{k-1} + ... + \binom{c_{1}}{1}$, for some $c_{k-1}>c_{k-2}>...>c_1$, then the maximum value for $a$ is between:
	
	\[\left[ \binom{c_{k-1}}{k} + ... + \binom{c_{1}}{2}\right] \leq f(r,k,b) \leq \left[ \binom{c_{k-1}}{k} + ... + \binom{c_{1}}{2}\right] + \tau.\]
\end{mythm}

\textbf{Remark:} Bollob\`{a}s and Eccless also proposed the weaker conjecture that $f\left(r,k,\binom{x}{k-1}\right)\leq \binom{x}{k}$ whenever $x$ is a positive real that makes $\binom{x}{k-1}$ an integer. We do end up proving this in the cases $k\leq 3$. $k=0,1,2$ are just corrolaries of theorem 1 while $k=3$ is found during the proof of the theorem 1. However, for the case $k\geq 4$, theorems 2 and 3 are still the best we have so far. 

\section{The case $k=0$}

This one is trivial and you'll obviously have $\mathcal{B}=\emptyset$ regardless of what $\mathcal{A}$ is.

\section{The case $k=1$}

This one is similarly trivial: the optimum will be $b=1$ regardless of what $a$ is. This can be achieved by letting $\mathcal{B}$ be an arbitrary set $B$ of size $r-1$, and $\mathcal{A}$ an arbitrary collection of $a$ sets of size $r$ all of which contain $B$.

\section{The case $k=2$}

For any pair of elements in $\mathcal{B}$, there will be at most 1 element in $\mathcal{A}$ that contains both (which if it exists will be their union). Since every element of $\mathcal{A}$ does contain a pair of elements of  $\mathcal{B}$, we have $|\mathcal{A}|\leq \binom{|\mathcal{B}|}{2}$.\\
\\
This is achieved by Bollobas's and Eccles's conjecture, ie, there is some $S$ of size $r-2$ (all of whose elements are larger than $b$), $\mathcal{B}$ is $\{S\cup\{i\} \, : \, i\leq b\}$ and $\mathcal{A}$ is $\{S\cup\{i,j\} \, : \, i,j\leq b\}$. Therefore $|\mathcal{A}|= \binom{|\mathcal{B}|}{2}$.

\section{The case $k=3$}

	Given a valid configuration ($\mathcal{A},\mathcal{B})$, we can construct a $k$-hypergraph with $b$ vertices corresponding to elements of $\mathcal{B}$. There are $a$ edges corresponding to elements of $\mathcal{A}$; each one of these contains at least $k$ elements of $\mathcal{B}$, and these $k$ vertices are going to be the vertices the edge is incident to (if there are more than $k$ of them, pick $k$ of them arbitrarily).\\
	
	Also, given two sets $B_1$ and $B_2$, we define the distance $d(B_1,B_2)$ between them as $|B_1\triangle B_2|/2$ (so two adjacent vertices are at distance 1 from each other). Note that given 2 such sets at distance 1, there is at most one single edge that contains both: $B_1 \cup B_2$.\\

\subsection{The case $b=\binom{c}{2}$}

Consider our $3$-hypergraph with $b$ vertices and $a$ edges. The average degree is $\frac{3a}{b}$. 

\paragraph{Paths of length 2}

Let the path of length 2, $P_2$ be the hypergraph consisting of 2 edges intersecting in a single point (which we'll call the \underline{center}), and with two distingued points, one on each edge that are not the center: $B_1$ and $B_2$; we'll call these the \underline{endpoints}. Given a copy of $P_2$ inside our graph, which we'll call $H$, let the center of $H$ be $c(H)$

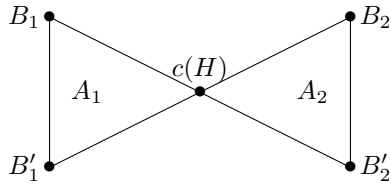
\begin{figure}[h]
\begin{tikzpicture}

\draw (0,0) node {$\bullet$};
\draw (-2,1) node {$\bullet$};
\draw (-2,-1) node {$\bullet$};
\draw (2,-1) node {$\bullet$};
\draw (2,1) node {$\bullet$};

\draw (0,0) node[above] {$c(H)$};
\draw (-2,1) node[left] {$B_1$};
\draw (-2,-1) node[left] {$B'_1$};
\draw (2,-1) node[right] {$B'_2$};
\draw (2,1) node[right] {$B_2$};

\draw (-1.5,0) node {$A_1$};
\draw (1.5,0) node {$A_2$};

\draw (0,0) -- (-2,1);
\draw (0,0) -- (-2,-1);
\draw (0,0) -- (2,-1);
\draw (0,0) -- (2,1);
\draw (-2,1) -- (-2,-1);
\draw (2,1) -- (2,-1);

\end{tikzpicture}
\centering
\caption{The hypergraph $P_2$}
\end{figure}

We want to count the number of $P_2$s. First we'll pick a pair of intersecting edges. The number of possible choices is:

\begin{eqnarray*}
& & |\{H : H \cong P_2\}|\\
&=& \sum_B |\{H  : H \cong P_2 \, , \, c(H)=B\}| \\
&=& \sum_B \binom{\deg(B)}{2}\\
&=&\sum_B \frac{\deg(B)^2}{2}- \frac{\deg(B)}{2}\\
& \geq & \frac{\left(\sum_B \deg(B)\right)^2}{2\sum_B 1}-\frac{\sum_B \deg(B)}{2}\\
&=&\frac{(3a)^2}{2b}-\frac{3a}{2}\\
&=&\frac{9a^2}{2b}-\frac{3a}{2}.\end{eqnarray*}

Set $\epsilon_1(B) = \left[\deg(B)-\sum_C \deg(C)/b\right]^2$. We have: 
\[\sum_{B\in\mathcal{B}} \epsilon_1(B) = \sum_{B\in\mathcal{B}} \deg(B)^2 -2\left[\sum_{B\in\mathcal{B}} \deg(B)\right]\left[\frac{\sum_{C\in\mathcal{B}} \deg(C)}{b}\right] + b\left[\frac{\sum_{C\in\mathcal{B}} \deg(C)}{b}\right]^2 = \sum_{B\in\mathcal{B}} \deg(B)^2 - \frac{\left(\sum_B \deg(B)\right)^2}{b} .\] 
This is exactly twice the error in the above inequality, so the number of intersecting pairs of edges is therefore exactly $\frac{9a^2/b-3a+\sum_{B\in\mathcal{B}} \epsilon_1(B)}{2}$.\\
\\
\paragraph{Paths of length 2 connecting points at distance 2} Now given a pair of intersecting edges, consider the number of pairs of points on it at distance 2 from each other. The maximum number of them is trivially 4 (2 choices on each side). It turns out that the minimum is in fact 2. Indeed, if $B_1$ and $B_2$ are on opposites sides of the path (see Figure 1) with $d(B_1,B_2)=1$, then we can write $B_1=c_H\cup\{x_1\}\backslash \{y_1\}$ and $B_2=c_H\cup\{x_2\}\backslash \{y_1\}$. Then $B'_2=c_H\cup\{x_2\}\backslash \{y_2\}$ so $B'_2$ and $B_1$ are at distance 2 and using a similar reasoning, $B_2$ and $B'_1$ are at distance 2. \\
\\
We are now interested in the number of paths $H$ of length 2 going from a vertex $B$ towards another vertex $B'$ at distance 2 from $B$ (note that the path from $B'$ to $B$ is counted as a seperate object from the path going in the opposite direction from $B$ to $B'$). We can call such an objects $J=\{(H,B_1,B_2): H \cong P_2 ; B_1,B_2\in H; d(B_1,B_2)=2\}$. How many such objects are there? Every path of length 2 either contributes 4,6 or 8 to the number. By default, we say each contributes 4, and there are a few special $P_2$s that contribute 2 or 4 extra.\\
\\
Let $\epsilon_2(B)$ be the number of pairs of intersecting edges containing $B$ such that $B$ is not the center and such that the number of pairs of points on it at distance 2 from each other is 3. Simiarly, let $\epsilon_3(B)$ be the number of copies of $P_2$ containing $B$ such that $B$ is not the center and such that the number of pairs of points on the pair of intersecting edges at distance 2 from each other is 4. Thus, $\sum_B \epsilon_2(B)/4$ is the number of pairs of intersecting edges that contribute 1 extra, while $\sum_B \epsilon_3(B)/4$ is the number of pairs of intersecting edges that contribute 2 extra. Thus:

\begin{eqnarray*}
|J| & = & 4|\{A_1,A_2\in\mathcal{A} : |A_1\cap A_2|=1\}| + \sum_{B\in\mathcal{B}} \epsilon_2(B)/2 + \sum_{B\in\mathcal{B}} \epsilon_3(B) \\
& = & 18a^2/b-6a+2\sum_{B\in\mathcal{B}} \epsilon_1(B) + \sum_{B\in\mathcal{B}} \epsilon_2(B)/2 + \sum_{B\in\mathcal{B}} \epsilon_3(B).
\end{eqnarray*}

\paragraph{Distance 2 pairs} Now given two points at distance 2, there are at most 4 paths of length 2 containing both (because there are at most 4 possible centers). Therefore the number of pairs of points $(B,B')$ is at least the number of paths from one vertex to another at distance 2 from it, divided by 4. \\
\\
Set $\epsilon_4(B)$ to be the number of points $B'$ at distance 2 from $B$ such that there are only 3 paths of length 2 between the two. $\epsilon_5(B)$ to be the number of points $B'$ at distance 2 from $B$ such that there are only 2 paths of length 2 between the two. $\epsilon_6(B)$ to be the number of points $B'$ at distance 2 from $B$ such that there is only 1 path fo length 2 between the two. $\epsilon_7(B)$ to be the number of points $B'$ at distance 2 from $B$ such that there are no paths of length 2 between the two. Using this, we get:

\[|J|  =   4|\{(B,B'):d(B,B')=2\}| - \sum_{B\in\mathcal{B}} \epsilon_4(B) - 2\sum_{B\in\mathcal{B}} \epsilon_5(B) - 3\sum_{B\in\mathcal{B}} \epsilon_6(B) - 4\sum_{B\in\mathcal{B}} \epsilon_5(B).\]

And thus the total number of ordered pairs of points $(B,B')$ at distance 2 from each other is:

\begin{eqnarray*}
|\{(B,B') : d(B,B')=2\}| & = &  \frac{|J| + \sum_{B\in\mathcal{B}} \epsilon_4(B) + 2\sum_{B\in\mathcal{B}} \epsilon_5(B) + 3\sum_{B\in\mathcal{B}} \epsilon_6(B) + 4\sum_{B\in\mathcal{B}} \epsilon_7(B)}{4}\\
 & = & \frac{9a^2}{2b} - \frac{3a}{2} + \sum_{B\in\mathcal{B}} \left[\frac{\epsilon_1(B)}{2} + \frac{\epsilon_2(B)}{8} + \frac{\epsilon_3(B)}{4} + \frac{\epsilon_4(B)}{4} + \frac{\epsilon_5(B)}{2} + \frac{3\epsilon_6(B)}{4} + \epsilon_7(B) .\right].
\end{eqnarray*}

We also know that the number of ordered pairs of points $(B,B')$ at distance 1 is at least $6a$ (6 from each edge). We set $\epsilon_8(B) $ to be the number of points at distance 1 from $B$ such that there is no edge connecting them. So $6a + \sum_B \epsilon_8(B)$ is the total number of ordered pairs at distance 1. Finally, let $\epsilon_9(B)$ be the number of points at distance at least 3 from $B$, so $\sum_B \epsilon_9(B)$ is the total number of ordered pairs of points at distance at least 3. Together with our calculated `number of ordered pairs at distance 2', this accounts for every possible ordered pair of points. We know that the total number of such pairs is $b(b-1)$. Thus:

\[b(b-1) = \frac{9a^2}{2b} - \frac{3a}{2} + 6a + \sum_{B\in\mathcal{B}} \left[\frac{\epsilon_1(B)}{2} + \frac{\epsilon_2(B)}{8} + \frac{\epsilon_3(B)}{4} + \frac{\epsilon_4(B)}{4} + \frac{\epsilon_5(B)}{2} + \frac{3\epsilon_6(B)}{4} + \epsilon_7(B) + \epsilon_8(B) + \epsilon_9(B)\right]. \]

We know that all the $\epsilon_i(B)$s are non-negative so we get the inequality: 

\[b(b-1) \geq \frac{9a^2/b-3a}{2}+6a.\]

Solving this for $a$ gives:

\begin{equation}
a\leq \frac{-3+\sqrt{8b+1}}{6}b . \label{poly}
\end{equation}

When we set $b=\binom{c}{2}$, this inequality gives us $a\leq \binom{c}{3}$. This is in fact tight. Indeed, we can look at the configuration from the hypothesis: let $S$ be a set of size $r-3$ (that doesn't contain any of 1,2,3,...,$c$) and let $\mathcal{B}=\{S\cup\{i,j\}|i,j\leq c\}$ and $\mathcal{A}=\{S\cup\{i,j,k\}|i,j,k\leq c\}$. Then this is a valid configuration and has $b=\binom{c}{2}$ and $a=\binom{c}{3}$.\\
\\
\textbf{Remark:} This formula is exactly the same as $f\left(r,3,\binom{x}{2}\right)\leq \binom{x}{3}$ for every real $x$ such that $\binom{x}{2}$ is a positive integer, thereby proving the weaker conjecture in the case $k=3$.

\subsection{The existance of a large nice hypergraph when $a$ is close to the upper bound}

For this section, we will assue that $c_2\geq 29$. If we set $b=\binom{c_2}{2}+c_1$ and $a=\binom{c_2}{3}+\binom{c_1}{2}+1$ for some $c_1<c_2$, then this means that the sum of all the $\epsilon_i(B)$ are small. In fact, we get: 

\[\frac{\sum_B}{b} \left[\frac{\epsilon_1(B)}{2} + \frac{\epsilon_2(B)}{8} + \frac{\epsilon_3(B)}{4} + \frac{\epsilon_4(B)}{4} + \frac{\epsilon_5(B)}{2} + \frac{3\epsilon_6(B)}{4} + \epsilon_7(B) + \epsilon_8(B) + \epsilon_9(B)\right]  = (b-1) - \frac{9a^2}{2b^2}-\frac{9a}{2b}. \]

This is an average, so in particular, there exists a vertex $B$ for which :

\[\frac{\epsilon_1(B)}{2} + \frac{\epsilon_2(B)}{8} + \frac{\epsilon_3(B)}{4} + \frac{\epsilon_4(B)}{4} + \frac{\epsilon_5(B)}{2} + \frac{3\epsilon_6(B)}{4} + \epsilon_7(B) + \epsilon_8(B) + \epsilon_9(B) \leq (b-1) - \frac{9a^2}{2b^2}-\frac{9a}{2b}. \]

And for brevity, we'll call the left hand side of this inequality $\gamma$. 

\begin{equation}
\gamma=\frac{\epsilon_1(B)}{2} + \frac{\epsilon_2(B)}{8} + \frac{\epsilon_3(B)}{4} + \frac{\epsilon_4(B)}{4} + \frac{\epsilon_5(B)}{2} + \frac{3\epsilon_6(B)}{4} + \epsilon_7(B) + \epsilon_8(B) + \epsilon_9(B) \label{defgamma}\\
\end{equation}
\[\gamma \leq \frac{c_2^2-c_2+2c_1-2}{2}-\frac{9}{2}\left(\frac{a}{b}\right)^2-\frac{9}{2}\left(\frac{a}{b}\right)\]

 We have $\frac{a}{b}=\frac{1}{3}\frac{c_2^3-3c_2^2+2c_2+3c_1^2-3c_1+6}{c_2^2-c_2+2c_1} = \frac{c_2-2-2c_1/c_2+3c_1^2/c_2^2-c_1/c_2^2+7c_1^2/c_2^3 -6c_1^3/c_2^4 }{3} + \frac{6-c_1/c_2+9c_1^2/c_2^2-20c_1^3/c_2^3+12c_1^4/c_2^4}{c_2^2-c_2+2c_1}$. But now the function $6-z+9z^2-20z^3+12z^4$ is always positive, so by removing the second term, we get $a/b \geq \frac{c_2-2-2c_1/c_2+3c_1^2/c_2^2-c_1/c_2^2+7c_1^2/c_2^3 -6c_1^3/c_2^4 }{3}$. The inequality now becomes: 

\begin{eqnarray*}
\gamma & \leq & \frac{c_2^2-c_2+2c_1-2}{2}-\frac{1}{2}\left(c_2-2-2c_1/c_2+3c_1^2/c_2^2-c_1/c_2^2+7c_1^2/c_2^3 -6c_1^3/c_2^4\right)^2 \\
& & -\frac{3}{2}\left(c_2-2-2c_1/c_2+3c_1^2/c_2^2-c_1/c_2^2+7c_1^2/c_2^3 -6c_1^3/c_2^4\right) \\
& = & \frac{c_2^2-c_2+2c_1-2}{2} - \frac{c_2^2-4c_2-4c_1+6c_1^2/c_2}{2}-\frac{4+6c_1/c_2+6c_1^2/c_2^2-24c_1^3/c_2^3+9c_1^4/c_2^4}{2} \\ 
& & -\frac{4c_1/c_2-24c_1^2/c_2^2-10c_1^3/c_2^3+66c_1^4/c_2^4-36c_1^5/c_2^5}{2c_2} -\frac{c_1^2/c_2^2  - 14 c_1^3/c_2^3 + 61 c_1^4/c_2^4 - 84 c_1^5/c_2^5+ 36 c_1^6/c_2^6}{2c_2^2} \\
& & -\frac{3c_2}{2}-\frac{-6-6c_1/c_2+9c_1^2/c_2^2}{2}-\frac{-3c_1/c_2+21c_1^2/c_2^2-18c_1^3/c_2^3}{2c_2} \\
& = & \frac{6c_1-6c_1^2/c_2}{2}+\frac{-15c_1^2/c_2^2+24c_1^3/c_2^3-9c_1^4/c_2^4}{2}+\frac{-c_1/c_2+3c_1^2/c_2^2+28c_1^3/c_2^3-66c_1^4/c_2^4+36c_1^5/c_2^5}{2c_2}\\
& & + \frac{-c_1^2/c_2^2  + 14 c_1^3/c_2^3 - 61 c_1^4/c_2^4 + 84 c_1^5/c_2^5- 36 c_1^6/c_2^6}{2c_2^2}
\end{eqnarray*}

Now since $c_1/c_2$ is between 0 and 1, we can look at the 3 last terms as functions of $c_1/c_2$ and see that each of of these functions is small on the interval $[0,1]$. To be more precise, the terms of order $O(1)$ are bounded above by $0.7041$, the terms of order $O(1/c_2)$ are bounded above by $0.39498/c_2$ and the terms of order $O(1/c_2^2)$ are bounded above by 0. Since $c_2\geq 24$, the sum of these terms is less than 1. So from this, we get that $\gamma$ is less than $\frac{6c_1-6c_1^2/c_2}{2}+1=3c_2.[c_1/c_2.(1-c_1/c_2)]+1$. Since the function $z(1-z)$ has a maximum of $1/4$, this is less than $3c_2/4+1$ and furthermore since $c_2\geq 29$, we get $\gamma\leq c_2$.  \\

\paragraph{Using $\gamma$ to restrict what the familes looks like}

We know from the definition of $\gamma$ (\ref{defgamma}) that $\epsilon_1(B) \leq 2\gamma$ so $|\deg(B)-3a/b|\leq \sqrt{2\gamma}$. So we know that the degree of $B$ is pretty close to its expected value.\\
\\
Now we also know that the number of points $B'$ at distance 2 from $B$ is $b-1 - 2\deg(B) - \epsilon_8(B)-\epsilon_9(B)$. We want to count the number of such points that have all 4 possible paths between it and $B$. This is just $b-1 -2\deg(B) - \epsilon_8(B)-\epsilon_9(B)-\epsilon_4(B)-\epsilon_5(B)-\epsilon_6(B)-\epsilon_7(B)$.\\

Again using the definition of $\gamma$ (\ref{defgamma}), we know that $\sum_{i=4}^9 \epsilon_i(B)$ is less than $4\gamma$ (with equality if and only if $\epsilon_4(B)=\gamma$ and all others are 0) so the number of points $B'$ at distance 2 from $B$ with all 4 possible paths between $B$ and $B'$ is at least $b-1-2\deg(B)-4\gamma$. We'll call such a configuration an \underline{octahedron}.\\
\\
We'll say that the $\deg(B)$ edges adjacent to $B$ are $B\cup\{x_1\}, B\cup\{x_2\},... B\cup\{x_{\deg(B)}\}$, and say that the $2\deg(B)$ vertices at distance 1 from $B$ on these edges are:\\ $B\cup\{x_1\}\backslash\{y_1\}, B\cup\{x_1\}\backslash\{z_1\}, B\cup\{x_2\}\backslash\{y_2\}, B\cup\{x_1\}\backslash\{z_2\},...,B\cup\{x_{\deg(B)}\}\backslash\{y_{\deg(B)}\}, B\cup\{x_{\deg(B)}\}\backslash\{z_{\deg(B)}\}$\\
\\
We'll now colour the edges incident to $B$ as follows: every edge incident to $B$ has vertices of the form $B$,$B\cup\{x\}\backslash\{y\}$ and $B\cup\{x\}\backslash\{z\}$. The colour of this edge is defined to be $\{y,z\}$. This colouring is useful because of its relationship to the octahedrons.

Indeed, an octahedron is formed of the 6 points: $B, B\cup\{x_i\}\backslash\{y_i\}, B\cup\{x_i\}\backslash\{z_i\}, B\cup\{x_j\}\backslash\{y_j\}, B\cup\{x_j\}\backslash\{z_j\}$ and $B\cup\{x_i,x_j\}\backslash\{y_i,z_i\}$ and also requires that $\{y_i,z_i\}=\{y_j,z_j\}$. That means that each octahedron contains 2 edges of the same colour. Furthermore, given a pair of edges of the same colour, there can only be at most 1 octahedron that contains both. So the number of octahedrons is smaller than the number of pairs of edges of the same colour. \\
\\
We'll define $s$ to be the size of the largest colour class. How large must $s$ be?

If $s$ is fixed and $s\geq \deg(B)/2$, then the maximum number of octahedrons is $\binom{s}{2}+\binom{\deg(B)-s}{2}=s^2+\frac{\deg(B)^2-\deg(B)}{2}-\deg(B)s$. We know that this quantity is less than $b-1-2\deg(B)-4\gamma$ so solving this equation in $s$ gives us:

\[s\geq \frac{\deg(B)+\sqrt{4b-4-16\gamma-6\deg(B)-\deg(B)^2}}{2}.\]

Using the fact that $\deg(B)\geq 3a/b-\sqrt{2\gamma}=\frac{c_2^3-3c_2^2+2c_2+3c_1^2-3c_1+6}{c_2^2-c_2+2c_1}-\sqrt{2\gamma} \geq c_2-3-\sqrt{2\gamma}$, we get that this is at least:

%\begin{eqnarray*}
% s & \geq & \frac{c_2-3-\sqrt{2\gamma}+\sqrt{c_2^2-2c_2+4c_1+5-18\gamma+c_2\sqrt{8\gamma}}}{2}\\
% & \geq & \frac{c_2-3-\sqrt{2\gamma}+c_2-1+2c_1/c_2-9\gamma/c_2+\sqrt{2\gamma}+o(1)}{2} \\
%& \geq & c_2-2+c_1/c_2-\frac{9}{2}\gamma/c_2+o(1)\\
%& \geq & c_2-2+c_1/c_2-\frac{27}{2}c_1/c_2(1-c_1/c_2)+o(1)\\
%&  \geq & c_2-9+o(1)
%\end{eqnarray*}

\begin{eqnarray*}
s & \geq & \frac{c_2-3-\sqrt{2\gamma}+\sqrt{c_2^2-2c_2+4c_1+5-18\gamma+c_2\sqrt{8\gamma}}}{2} \\
 & \geq & \frac{c_2-3-\sqrt{2\gamma}+\sqrt{c_2^2-2c_2+5-18\gamma+c_2\sqrt{8\gamma}}}{2}.
\end{eqnarray*}

Now we split into 2 cases depending on if $\gamma< 8$ or $\gamma\geq 8$. If $\gamma<8$, then:

\begin{eqnarray*}
s & > & \frac{c_2-3-\sqrt{2.8}+\sqrt{c_2^2-22c_2+5-18.8}}{2}\\
 & = & \frac{c_2-7+\sqrt{c_2^2-2c_2-139}}{2}\\
  & \geq & \frac{c_2-7+\sqrt{c_2^2-7c_2+14}}{2}\\
  & \geq & \frac{c_2-7+\sqrt{c_2^2-7c_2+49/4}}{2}\\
  & \geq & \frac{c_2-7+|c_2-7/2|}{2}\\
  & = & \frac{2c_2-21/2}{2}\\
  & = & c_2-21/4.\\
\end{eqnarray*}

If on the other hand we have $\gamma\geq 8$, then $\sqrt{2\gamma}\geq 4$ and therefore:

\begin{eqnarray*}
s & \geq & \frac{c_2-3-\sqrt{2\gamma}+\sqrt{c_2^2-22c_2+5+2\gamma+c_2\sqrt{8\gamma}}}{2}\\
 & \geq & \frac{c_2-3-\sqrt{2\gamma}+\sqrt{c_2^2-22c_2+5+2\gamma+c_2\sqrt{8\gamma}-2(c_2-29)-24(\sqrt{2\gamma}-4)}}{2}\\
 & = & \frac{c_2-3-\sqrt{2\gamma}+\sqrt{c_2^2-24c_2+159+2\gamma+c_2\sqrt{8\gamma}-24\sqrt{2\gamma}}}{2}\\ 
  & \geq & \frac{c_2-3-\sqrt{2\gamma}+\sqrt{c_2^2-24c_2+144+2\gamma+c_2\sqrt{8\gamma}-24\sqrt{2\gamma}}}{2}\\
 & = & \frac{c_2-3-\sqrt{2\gamma}+|c_2-12+\sqrt{2\gamma}|}{2}\\
 & = & \frac{2c_2-15}{2}\\
  & = & c_2-15/2.\\
\end{eqnarray*}

Since $s$ and $c_2$ are both integers, we actually get $s\geq c_2-7$, and this holds regardless of which case we were in. So this colour class of size $s$ encompases most of the neighbourhood of $B$ when $c_2\geq 24$. In fact, there are only at most $\deg(B)+7-c_2$ points of the neighbourhood that are outside. Since $\deg(B)\leq \frac{3a}{b}+\sqrt{2\gamma}\leq \frac{3\binom{c_2+1}{3}}{\binom{c_2}{2}}+\sqrt{2\gamma}=c_2+1+\sqrt{2\gamma}$, that means there are at most  $8+\sqrt{2\gamma}$ points in the neighbourhood of $B$ outside our colour class. We'll say that the colour of this large colour class is $\{y,z\}$. Finally, we'll define $S=B\backslash\{y,z\}$ and define the \underline{nice hypergraph} to be the set of all vertices an edges that contain $S$ as a subset. Our nice hypergraph contains our large colour class as well as all the octahedrons related to it.\\
\\
 We ask ourselves how many octahedrons do we have in the nice hypergraph? We know that there are at least $b-1-2\deg(B)-4\gamma$ octahedrons that contain $B$ in total. The maximum number of octahedrons containing $B$ that are not in the nice hypergraph is $\binom{\deg(B)-s}{2} \leq \binom{8+\sqrt{2\gamma}}{2} = 28+17\sqrt{\gamma/2}+\gamma$. Therefore there are at least $b-1-2\deg(B)-4\gamma-(28+17\sqrt{\gamma/2}+\gamma)\geq b-29-2\deg(B)-5\gamma-17\sqrt{\gamma/2}$ octahedrons in our big colour class. That's almost all the vertices in the graph! Remembering the $2\deg(B)$ vertices adjacent to $B$ and $B$ itself, that means there are only $28+17\sqrt{\gamma/2}+5\gamma < 9c_2$ left unaccounted for. Before we finish up the proof, we will do the following easy example to demonstrate what we have so far: \\

\subsection{Example case where $b=\binom{c_2}{2}$ and $a=\binom{c_2}{3}$}

First of all, we check that $\gamma=0$. This implies that the degree of every vertex is exactly $3a/b=c_2-2$ and for every vertex $B$, the number of vertices at distance 2 from it is exactly $b-1-2\deg(B)$. Finally, we know that the every single one of these points forms an octahedron with $B$ and 2 points in the neighbourhood of $B$.   \\
\\
We can write the neighbourhood of $B$ as $B\cup\{x_i\}\backslash\{y\}$ and $B\cup\{x_i\}\backslash\{z\}$ for all $i$ between $1$ and $c_2-2$. Finally, we can write every other point of the graph in the form $B\cup\{x_i,x_j\}\backslash\{y,z\}$ for all $1\leq i<j\leq c_2-2$. \\
\\
Let $S=B\cup\backslash\{y\}$ and then $\mathcal{A}$ is exactly the family $\{S\cup \{t_1,t_2,t_3\} | t_1,t_2,t_3\in\{y,z,x_1,x_2,...,x_{c_2-2}\}\}$ while $\mathcal{B}$ is exactly the family $\{S\cup \{t_1,t_2\} | t_1,t_2\in\{y,z,x_1,x_2,...,x_{c_2-2}\}\}$.

\begin{mylem}
The optimal configuration when $k=3$, $b=\binom{c_2}{2}$ and $a=\binom{c}{3}$ is of the form:\\
$\mathcal{A} = \{S\cup T | T\in\{x_1,x_2,...,x_{c_2}\}^{(3)}\}$\\
$\mathcal{B} = \{S\cup T | T\in\{x_1,x_2,...,x_{c_2}\}^{(2)}\}$.
\end{mylem}

\subsection{Other cases}
 
As a reminder, at this stage we know that there is a large `nice hypergraph' of vertices that all contain $S$ as subset. We'll say that there are exactly $\beta$ vertices not in our nice hypergraph (which leaves $b-\beta$ vertices in the nice hypergraph). We know that $\beta < 9c_2$. How many edges can there be in this graph now that we have this information? There are 3 types of edges, depending on how many vertices are in the nice hypergraph:\\
\\
$\bullet$ The edges that are entirely contained in the nice hypergraph. We shall call these \underline{nice edges}. All the vertices in the nice hypergraph contain $S$ so by the Kruskal Katona Theorem, the most edges entirely contained within is when they form an initial segment of the colex ordering. So if we write $b-\beta = \binom{d_2}{2}+\binom{d_1}{1}$, then we have at most $\binom{d_2}{3}+\binom{d_1}{2}$ nice edges. Because $\beta<9c_2$ and $c_2\geq 24$, we have $\binom{d_2+1}{2} > b-\beta > \binom{c_2}{2}-9c_2=\frac{c_2^2-19c_2}{2}\geq \frac{c_2^2-27c_2+192}{2} > \binom{c_2-13}{2}$. Therefore $d_2+1>c_2-13$ and because these are integers, $d_2\geq c_2-13$.\\
$\bullet$ The edges that that contain 1 or 2 vertices from the nice hypergraph and 2 or 1 from outside. We shall call these \underline{linking edges}. There are at most $\beta$ of them because given any point $T$ outside the nice hypergraph, the only potential edge that can connect to elements in the nice hypergraph is $T\cup S$.\\
$\bullet$ The edges entirely outside the nice hypergraph. We'll call these \underline{outside edges}. If we just apply our earlier result (\ref{poly}), we get that there are at most $\frac{-3+\sqrt{8\beta+1}}{6}\beta$ of them.\\
 \\
In total, we have at most $\binom{d_2}{3}+\binom{d_1}{2}+\frac{3+\sqrt{8\beta+1}}{6}\beta$ edges in our hypergraph. This is also equal to $\binom{c_2}{3}+\binom{c_1}{2}+1$, so we get:

\[(c_2-d_2)\frac{c_2^2+c_2d_2+d_2^2-3c_2-3d_2+2}{6} + (c_1-d_1)\frac{c_1+d_1-1}{2} \leq \frac{3+\sqrt{8\beta+1}}{6}\beta.\]

And $\beta = (c_2-d_2)\frac{c_2+d_2-1}{2} + (c_1-d_1)$. Let $\phi=c_2-d_2$; we know that $0\leq \phi\leq 13$. Replace all instances of $d_2$ with $c_2-\phi$ in the inequality. We end up with: 

\begin{equation}
\phi\left(3c_2^2-3\phi c_2+\phi^2-6c_2+3\phi+2\right) + 3(c_1-d_1)(c_1+d_1-1) \leq (3+\sqrt{8\beta+1})\beta. \label{inequality} \end{equation}

And $\beta = \phi\left(c_2-\frac{\phi+1}{2}\right) + (c_1-d_1)$.

\paragraph{Case 1: $2\leq \phi\leq 13$\\}

We know that $0\leq c_1\leq c_2-1$ and $0\leq d_1\leq d_2-1$ so $(c_1-d_1)(c_1+d_1-1)\geq -(d_2-1)(d_2-2)=-(c_2-1)(c_2-2)+\phi(2c_2-3)-\phi^2$ and $\beta\leq \phi\left(c_2-\frac{\phi+1}{2}\right) + c_2-1=(\phi+1)(c_2-\phi/2)-1$. Putting these back into inequality \ref{inequality}, we get:\\

\begin{eqnarray*}
& & \phi\left(3c_2^2-3\phi c_2+\phi^2-6c_2+3\phi+2\right) -3(c_2-1)(c_2-2)+3\phi(2c_2-3)-3\phi^2 \\ 
& \leq & (3+\sqrt{8(\phi+1)(c_2-\phi/2)-7})[(\phi+1)(c_2-\phi/2)-1].
\end{eqnarray*}
So:
\begin{eqnarray*}
& & (3\phi-3)c_2^2+(-3\phi^2+9)c_2+(\phi^3-7\phi-6) \\
& \leq & 3(\phi+1)(c_2-\phi/2)-3+\sqrt{8(\phi+1)(c_2-\phi/2)-7}[(\phi+1)(c_2-\phi/2)-1]. \\
\end{eqnarray*}
Therefore:

\begin{eqnarray*}
& & (3\phi-3)c_2^2+(-3\phi^2-3\phi+6)c_2+(\phi^3+3/2\phi^2-11/2\phi-3) \\
& \leq & \sqrt{8(\phi+1)(c_2-\phi/2)-7}[(\phi+1)(c_2-\phi/2)-1]. \\
\end{eqnarray*}

If we examine this function of $c_2$ for $\phi$ between 2 and 13, we find that this only holds when $c_2 < 29$. We assumed that $c_2\geq 29$ so therefore this case cannot occur.

\paragraph{Case 2: $\phi=1$\\} Then inequality \ref{inequality} becomes:

\begin{equation}
\left(3c_2^2-9c_2+6\right) + (c_1-d_1)\left(3c_1+3d_1-3\right) \leq  (3+\sqrt{8\beta+1})\beta. \label{case2}
\end{equation}

where $\beta = c_2-1+(c_1-d_1) \leq 2c_2-2$ so that implies: 

\[\left(3c_2^2-12c_2+9\right) + 3c_1^2 -3d_1^2 +6(d_1-c_1) \leq  8c_2^{3/2} . \]

But now $3c_2^2-8c_2^{3/2}-12c_2+9\geq 3(c_2-4/3\sqrt{c_2}-4)^2$ so this implies that $d_1\geq c_2-4/3\sqrt{c_2}-4$. We also get that $c_1\leq \frac{1}{3}\sqrt{8c_2^{3/2}+12c_2-9} < c_2/2$.

%When $c_2$ is large, the only way that this inequality is satisfied is if $3c_2^2+3c_1^2-3d_1^2=O(c_2\sqrt{\beta})$ so $d_1=c_2-O(\sqrt{\beta})$ and $c_1=O(c_2^{1/2}\sqrt[4]{\beta})$. \\
Replacing $d_1=c_2-\epsilon$ in inequality \ref{case2}, we get: 

\[(6\epsilon-6)c_2+3c_1^2-6c_1-3\epsilon^2-6\epsilon+9 \leq \sqrt{8c_1+8\epsilon-7}(c_1+\epsilon-1) .\]

And we know that $2\leq \epsilon\leq 4/3\sqrt{c_2}+4< c_2/2$. Since $\sqrt{8c_1+8\epsilon-7}<\sqrt{8c_2}$, we can say that:

\[-3\epsilon^2+\epsilon(6c_2-\sqrt{8c_2}-6) +3(c_1-1)^2-\sqrt{8c_2}(c_1-1)-6c_2+6 \leq 0 . \]
\\
The value for $c_1$ that minimises the left hand side is $c_1=1+\sqrt{2c_2}/3$ so we can without loss of generality assume that is what $c_1$ is, and that gives us:

\[-3\epsilon^2+\epsilon(6c_2-\sqrt{8c_2}-6) -20/3c_2+6 \leq 0 . \]

If we set $\epsilon$ to be at its maximum value: $\epsilon=4/3\sqrt{c_2}+4$, we get the left hand side is $8c_2^{3/2}+(12-8\sqrt{2}/3)c_2-(40+8\sqrt{2})\sqrt{c_2}-2$ which is positive for $c_2\geq 24$. Similarly, when $\epsilon$ is at its minimium value, $\epsilon=2$, we have $-18+16/3c_2-4\sqrt(2c_2)$ which is also positive for $c_2\geq 24$. As this function is a quadratic polynomial with a negative leading term, it is concave and therefore this inequality does not hold for any valid $\epsilon$. So therefore this case cannot occur.

\paragraph{Case 3: $c_2=d_2$ but $c_1\neq d_1$\\} Then we can simplify (\ref{inequality}) further to $(3c_1+3d_1-6) \leq \sqrt{8c_1-8d_1+1}$ so either $c_1+d_1\leq 2$ or $3c_1-6\leq \sqrt{8c_1+1}$ so $9c_1^2-43c_1+35 \leq 0$ so $c_1\leq 3$ in either case. Looking at each subcase individually, we get the following cases:\\
\\
$\bullet$ Case 3.1. $d_1=0$ and $c_1 = 3$. In this case, our potential counter-example consists of 3 vertices outside the nice hypergraph with one outside edge and three linking edges (one per vertex), which is one more than the $\binom{3}{2}=3$ we expected. But this configuration is actually impossible.\\
Indeed, if we remove these 3 outside vertices and their 4 edges, we're left with a $\binom{c_2}{2}$ vertices and $\binom{c_2}{3}$ edges so by lemma 1, there is only a single unique solution: $\mathcal{B} = \{S\cup T | T\in\{x_1,x_2,...,x_{c_2}\}^{(2)}\}$ and $\mathcal{A} = \{S\cup T | T\in\{x_1,x_2,...,x_{c_2}\}^{(3)}\}$.\\
Also note that each of the 3 linking edges is only incident to one of the outside vertices, which means that its two other endpoints are inside the nice hypergraph. Say they are $S\cup\{x_1,x_2\}$ and $S\cup\{x_1,x_3\}$. Then the linking edge has to be $S\cup\{x_1,x_2,x_3\}$. But this is one of the nice edges that we've already counted. Contradiction. Therefore this configuration is indeed impossible.\\
\\
$\bullet$ Case 3.2. $d_1=0$ and $c_1 = 2$. In this case, our potential counter example consists of 2 vertices outside the nice hypergraph with two linking edges, which is one more than the $\binom{2}{2}=1$ we expected. This configuration is also impossible. \\
Similarly to the last bullet point, we note that if we remove the 2 extra vertices and edges, we are left with families of the type: $\mathcal{B} = \{S\cup T | T\in\{x_1,x_2,...,x_{c_2}\}^{(2)}\}$ and $\mathcal{A} = \{S\cup T | T\in\{x_1,x_2,...,x_{c_2}\}^{(3)}\}$. But also each of the 2 linking edges has to be of the form $S\cup\{x_1,x_2,x_3\}$, which is not a linking edge at all, but rather a nice edge that we have already counted. Contradiction.\\
\\
$\bullet$ Case 3.3. $d_1=0$ and $c_1 = 1$. In this case, our potential counter example consists of 1 vertex outside the nice hypergraph with one linking edges, which is one more than the $\binom{1}{2}=0$ we expected. This configuration is also impossible and the proof is identical to the last two bullet points.\\
\\
$\bullet$ Case 3.4. $d_1=1$ and $c_1=2$. In this case, our potential counterexample can actually work. However it has $\binom{c_2}{3}+\binom{1}{2}+1$ edges which is the same as $\binom{c_2}{3}+\binom{2}{2}$ that our usual example gives us so it doesn't give any improvement.
\\
\\
\\
\\
So in conclusion, we have proved the conjecture in the case $k=3$ and for $c_2\geq 29$. This finishes the proof of Theorem 1.\\
\\
\\Remark: we did not answer the question of what happens when $c_2<29$; however, this is only finitely many cases so it could in theory be solved by simply checking all the cases individually.

\section{The case $k\geq 4$}

Similarly to what we did in the case $k=3$, we will define a $k$-hypergraph, whose vertices are the elements of $\mathcal{B}$, and whose edges are the elements of $\mathcal{A}$. We know that every element $A$ of $\mathcal{A}$ contains at least $k$ elements of $\mathcal{B}$ as subsets, so pick $k$ of them arbitrarily and they will form the vertices incident to $A$. When a vertex $B$ is incident to an edge $A$, we will say that $B\prec A$. This relation is almost the same as the subset relation $B\subset A$ ; the only differences is that we have excluded a few of them arbitrarily to make every $A\in\mathcal{A}$ be related to exactly $k$ elements of $\mathcal{B}$.\\
\\
Similarly to the case $k=3$, we again define the distance $d(B_1,B_2)=|B_1\triangle B_2|/2$. \\
\\
The proof will follow what we did in the case $k=3$ except we will be looking at pairs of points at distance $k-1$ from each other, instead of pairs at distance 2, and the paths joining them together will have length $k-1$ instead of length 2. 

\begin{mydef}

A \underline{path of length $i$}, $P_i$, is a sequence of alternating vertices and edges in a hypergraph such that every vertex is on the previous edge, and every edge contains the previous vertex and where there are $i+1$ vertices and $i$ edges. The vertices and edges are not necessarily distinct.\\
\\
The first vertex and the last vertex of a path are called the \underline{endpoints}.

\end{mydef}

\noindent \textbf{Note:} The definition of $P_2$ is slightly different from the one we used in the previous section as we now allow self-intersections.\\
\\
There are also a few complications in how to count these paths when $i\geq 3$. We shall use a method developed by Szegedy \cite{Szegedy} that is used to solve the Sidorenko conjecture \cite{Sidorenko} in certain special cases. The Sidorenko conjecture states that if $H$ is a bipartite graph with $e(H)$ edges and $G$ is a graph with $n$ vertices and average degree $d$, then the number of homomorphisms from $H$ to $G$ is at least $nd^{e(H)}$. The conjecture has not been solved completely, but Szegedy's work proves it in a lot of cases and will be sufficient for our purposes. We shall first adapt Szegedy's work to count the number of paths of length $i$, then we will make a few more modifications that allow us to count the number of paths of length $i$ whose endpoints are at distance $i$ from each other.

\subsection{We define a probability distribution}

In our $k$-hypergraph, let $L_i$ be the set of paths of length $i$. For simplicity, the path going from $B$ to $B'$ and the opposite path going from $B'$ to $B$ are not considered the same path. Let $M_i$ be the subset of $L_i$ that contains only those paths whose end points are at distance $i$ from each other (so notably do not have any self-intersections)\\
\\
We are going to define a probability distribution $\mu_i: L_i \rightarrow [0,1]$ on each these objects, by induction. An element of $L_1$ is just a single edge with two distingued vertices on it. There are $a$ edges and then $k^2$ choices for $B$ and $B'$, so there are $ak^2$ elements of $L_1$. $\mu_1$ is going to pick one uniformly at random, so $\mu_1(l)=\frac{1}{ak^2}$ for any $l\in L_1$.\\
\\
For $i>1$, pick $l\in L_i$. This is a path from $B$ to $B'$. Remove $B'$ and its incident edge (which we'll call $A$) to get a shorter path $l'\in L_{i-1}$ going from $B$ to some other vertex $B''$. We define $\mu_i(l)=\frac{\mu_{i-1}(l')}{\deg(B'') \cdot k}$. \\
\\
Alternatively, if we write the vertices of the path as $B_0$, $B_1$, ..., $B_i$, then $\mu_i(l)=\frac{1}{\deg(B_1)\deg(B_2)...\deg(B_{i-1})ak^{i+1}}$.\\
\\
For a vertex $B$, we will also define $\mu_0(B)=\frac{\deg(B)}{ka}$. This is a probability distribution and moreover, if you consider a single vertex to be a path of length 0, then it agrees with all the other $\mu_i$s and satisfies the same properties. \\
\\
\\
\subsection{Size of $L_i$}
We claim that this probability distribution has the following property for any given vertex $B'$: 
\begin{equation}
\sum_{l\in L_i \text{ ending at $B'$}} \mu_i(l) = \mu_0(B) \label{prop1}
\end{equation} Indeed, when $i=1$, we have:

\[\sum_{l\in L_1 \text{ ending at $B'$}} \mu_1(l) = \sum_{A \succ B'; B\prec A} \frac{1}{ak^2} = \sum_{A \succ B'} \frac{1}{ak} = \frac{\deg(B)}{ka} = \mu_0(B).\]

For larger $i$, we have:

\[ \sum_{l\in L_i \text{ ending at $B'$}} \mu_i(l) = \sum_{\substack{l'\in L_{i-1} \text{ending in B''} \\ A\succ B'\, ; B''\prec A }} \frac{\mu_{i-1}(l')}{\deg(B'') \cdot k}.\]

Using the induction hypothesis (\ref{prop1}), we get that this is equal to:

\[\sum_{A\succ B'\, ; B'' \prec A } \frac{\deg(B'')}{ak}\frac{1}{\deg(B'') \cdot k} = \sum_{A\succ B'} \frac{1}{ak} = \mu_0(B).\]

So by induction, the claim is proved.\\
\\
\\
\\
We'll now define $D(\mu_i)=\sum_{l\in L_i}-\ln(\mu_i(l))\mu_i(l)$. Note that this is an entropy, so notably it is maximal when $\mu_i$ is uniform on $L_i$. Therefore $D(\mu_i)\leq \sum_{l\in L_i} -\ln(\frac{1}{|L_i|})\frac{1}{|L_i|} = \ln(|L_i|)$. \\
\\
We now want to calculate $D(\mu_i)$ to get a lower bound on $|L_i|$. For $i=1$, we have $\mu_1(l)=\frac{1}{ak^2}$ for all $l$, so $D(\mu_1)=\sum_{l\in L_1}-\ln(\frac{1}{ak^2})\frac{1}{ak^2}=\ln(ak^2)$.\\
\\
For larger $i$, we have, for every element $l\in L_i$, $l=l'\cup\{A,B'\}$, where $l'$ is the path of length $i-1$ consisting of the first $i$ vertices and $i-1$ edges, $A$ is the additional edge and $B'$ is the additional vertex. We'll say that $A$ connects to $l$ at the vertex $B''$.  Using the definition of $\mu_i$, we can rewrite $D(\mu_i)$ as:

\begin{eqnarray*} 
D(\mu_i)& = & \sum_{\substack{l'\in L_{i-1} \text{ ending at $B''$} \\ A\succ B'' \\ B'\prec A}}-\ln\left(\frac{\mu_{i-1}(l')}{\deg(B'')\cdot k}\right)\frac{\mu_{i-1}(l')}{\deg(B'')\cdot k}\\ 
&=&\sum_{\substack{l'\in L_{i-1} \text{ ending at $B''$}\\ A\succ B''}}-\ln\left(\frac{\mu_{i-1}(l')}{\deg(B'')\cdot k}\right)\frac{\mu_{i-1}(l')}{\deg(B'')} \\
&=&\sum_{\substack{l'\in L_{i-1} \text{ ending at $B''$}\\ A\succ B''}}\left[-\ln\left(\mu_{i-1}(l')\right)\frac{\mu_{i-1}(l')}{\deg(B'')}\right]+\left[\ln\left(\deg(B'')\right)\frac{\mu_{i-1}(l')}{\deg(B'')}\right]+\left[\ln\left(k\right)\frac{\mu_{i-1}(l')}{\deg(B'')}\right]\\
\end{eqnarray*}
Using the definition of $D$, property (\ref{prop1}), and the fact that $\mu_{i-1}$ is a probability distribution hence sums to 1, we can rewrite these 3 terms as:
\begin{eqnarray*}
D(\mu_i)&=&\left[\sum_{l'\in L_{i-1} \text{ ending at $B''$}}-\ln\left(\mu_{i-1}(l')\right)\mu_{i-1}(l')\right]+\left[\sum_{B'' ; A\succ B''}\ln(\deg(B''))\frac{\deg(B'')/ak}{\deg(B'')}\right] + \left[\ln(k)\right]\\
&=&\left[D(\mu_{i-1})\right]+\left[\sum_{B''}\ln\left(\frac{\deg(B'')}{ak}\right)\frac{\deg(B')}{ak}+\sum_{B''}\ln(ak)\frac{\deg(B'')}{ak}\right] + \left[\ln(k)\right]\\
&= & D(\mu_{i-1}) - D(\mu_0) +\ln(ak)+\ln(k)\\
&=&D(\mu_{i-1})- D(\mu_0) +\ln(ak^2).
\end{eqnarray*}

Now we use the entropy inequality on $\mu_0$ to say that $D(\mu_0)\leq \ln(b)$ and this gives us:

\[D(\mu_i) \geq D(\mu_{i-1})-\ln(ak^2/b).\] 

So by induction, we have $D(\mu_{i})\geq\ln(ak^2)+(i-1)\ln(ak^2/b) = \ln(\left(\frac{ak^2}{b}\right)^i b)$ and therefore we have at least $\left(\frac{ak^2}{b}\right)^i b$ paths of length $i$.\\
\\
\\
\subsection{Probability of being in $M_i$}

But we want to know the size of $M_i$, not $L_i$, so we need to do more work. Given a random element $l$ of $L_i$ (chosen according to the probability distribution $\mu_i$), what is the probability that it is actually in $M_i$? We'll call this probability $\mathbb{P}_i$. We know that $l$ is in $M_i$ if and only if its vertices are of the form $B$, $B\cup\{x_1\}\backslash\{y_1\}$, $B\cup\{x_1,x_2\}\backslash\{y_1,y_2\}$, ..., $B\cup\{x_1,x_2,...,x_i\}\backslash\{y_1,y_2,...,y_i\}$ while its edges are of the form $B\cup\{x_1\}$, $B\cup\{x_1,x_2\}\backslash\{y_1\}$, $B\cup\{x_1,x_2,x_3\}\backslash\{y_1,y_2\}$, ..., $B\cup\{x_1,x_2,...,x_i\}\backslash\{y_1,y_2,...,y_{i-1}\}$. \\

The first $i$ vertices and $i-1$ edges are of the correct form with probability $\mathbb{P}_{i-1}$, so we only need to consider the last vertex and last edge. The last edge has to contain $B''=B\cup\{x_1,x_2,...,x_{i-1}\}\backslash\{y_1,y_2,...,y_{i-1}\}$. There are only $i-1$ ways this fails to be of the correct form: $A=B''\cup\{y_1\}$, $A=B''\cup\{y_2\}$, ..., $A=B''\cup\{y_{i-1}\}$. The probability of this occuring is at most $\mu_{i-1}(l')\frac{i-1}{\deg(B'')}$.\\
Now suppose the last edge $A$ is of the correct form. How many ways can we place the $B'$ incorrectly? Well there are only $i$ possibilities: $A\backslash\{x_1\}$,$A\backslash\{x_2\}$, ..., $A\backslash\{x_i\}$. The probability this occurs is at most $\mu_{i-1}(l')\frac{\deg(B'')-i+1}{\deg(B'')}\frac{i}{k}$. So the total probability of success is $\mu_{i-1}(l')\frac{\deg(B'')-i+1}{\deg(B'')}\frac{k-i}{k}$.\\
 In fact, without loss of generality, we can actually assume that there are always exactly $i$ vertices that fail the second step because if there are fewer, just pick some of them arbitrarily and declare them failures. Thus, we can assume without loss of generality that $\mathbb{P}_i\leq \frac{k-1}{k}\mathbb{P}_{i-1}$ and hence $\mathbb{P}_i\leq \frac{(k-1)!}{(k-i-1)!k^i}$.  \\

So the total probability that $l$ is in $M_i$ is at least : 

\begin{eqnarray*}
& & \sum_{l'\in M_{i-1} }\mu_{i-1}(l')\frac{(\deg(B'')-i+1)(k-i)}{\deg(B'')k} \\
&=& \frac{k-i}{k}\mathbb{P}_{i-1}-\sum_{l'\in M_{i-1} }\mu_{i-1}(l')\frac{(i-1)(k-i)}{\deg(B'')k} \\
&\geq & \frac{k-i}{k}\mathbb{P}_{i-1}-\sum_{l'\in L_{i-1} }\mu_{i-1}(l')\frac{(i-1)(k-i)}{\deg(B'')k}\\
&= & \frac{k-i}{k}\mathbb{P}_{i-1}-\sum_{B''}\frac{(i-1)(k-i)}{ka.k} \\
& =&\frac{k-i}{k}\mathbb{P}_{i-1}-\frac{(i-1)(k-i)}{k}\frac{b}{ka}.
\end{eqnarray*}

(We use property (\ref{prop1}) again to go from the third line to the fourth.) By induction, we claim that $\mathbb{P}_i\geq \frac{(k-1)!}{(k-i-1)!k^i} - \left[\frac{k-i}{k}.\frac{i(i-1)}{2}\right]\frac{b}{ka}$. For $i=1$, we have $P_1=\frac{k-1}{k}$, which agrees with the formula. For larger $i$, we have from the above inequality:

\begin{eqnarray*}
\mathbb{P}_i & \geq & \frac{k-i}{k}\frac{(k-1)!}{(k-i)!k^{i-1}}-\frac{k-i}{k}\left[\frac{k-i+1}{k}.\frac{(i-1)(i-2)}{2}\right]\frac{b}{ka}-\frac{(i-1)(k-i)}{k}\frac{b}{ka}\\
& \geq & \frac{(k-1)!}{(k-i-1)!k^i} - \frac{k-i}{k}\left[\frac{(i-1)(i-2)}{2}\right]\frac{b}{ka}-\frac{k-i}{k}[i-1]\frac{b}{ka}\\
&=&\frac{(k-1)!}{(k-i-1)!k^i} - \left[\frac{k-i}{k}.\frac{i(i-1)}{2}\right]\frac{b}{ka}.
\end{eqnarray*}

So by induction, we know that: \begin{equation}
\mathbb{P}_i\geq \frac{(k-1)!}{(k-i-1)!k^i} - \left[\frac{k-i}{k}.\frac{i(i-1)}{2}\right]\frac{b}{ka}. \label{prob}
\end{equation}
\\
\\
\subsection{Size of $M_i$}

\paragraph{Step 1: Setting up a proof by induction\\}

So we know the number of paths, and we know the probability that one of our paths is `straight'. From this, we want to find the number of paths that are actually straight. This is slightly more complicated than it seems because the probability is not uniform. However, we can still find a lower bound. First, we do another entropy inequality:

\begin{eqnarray*}
\ln(|M_i|) &=& \sum_{l\in M_i} -\ln\left(\frac{1}{|M_i|}\right)\frac{1}{|M_i|} \\
& \geq & \sum_{l\in M_i} -\ln\left(\frac{\mu(l)}{\mathbb{P}_{i}}\right)\frac{\mu_i(l)}{\mathbb{P}_{i}}\\
&=& \ln(\mathbb{P}_i)+\frac{\sum_{l\in M_i} -\ln(\mu_i(l))\mu_i(l)}{\mathbb{P}_i}.
\end{eqnarray*}

So now we want to know find a lower bound on $\sum_{l\in M_i} -\ln(\mu_i(l))\mu_i(l)$. We claim the following: 
\begin{equation}\sum_{l\in M_i} -\ln(\mu_i(l))\mu_i(l) \geq \frac{(k-1)!}{(k-i-1)!k^i}\ln\left(\left(\frac{ak^2}{b}\right)^ib\right) - O\left(\frac{b\ln(b)}{ka}\right). \label{prop2}
\end{equation}
If that ends up being true, that would imply that

\begin{equation}
\ln(|M_i|) \geq \ln(\mathbb{P}_i)+\frac{(k-1)!}{(k-i-1)!k^i\mathbb{P}_i}\ln\left(\frac{ak^2}{b}\right) - O\left(\frac{b\ln(b)}{ka}\right). \label{prop3}
\end{equation}

We will prove (\ref{prop2}) by induction on $i$. For $i=1$, we get $\sum_{l\in M_1} -\ln(\mu_1(l))\mu_1(l) = \sum_{l\in M_1} -\ln\left(\frac{1}{ak^2}\right)\frac{1}{ak^2} = \ln(ak^2)\frac{a(k-1)k}{ak^2} = \frac{(k-1)!}{(k-2)!k^1}\ln\left(\left(\frac{ak^2}{b}\right)^1b\right)$ so it is true for $i=1$.

\paragraph{Step 2: For larger $i$, expressing $\sum_{l\in M_i}-\ln(\mu_i(l))\mu_i(l)$ as a sum of five terms\\}

We proceed as follows:

\begin{eqnarray*}
& & \sum_{l\in M_i} -\ln(\mu_i(l))\mu_i(l)\\
&=& \sum_{\substack{B''\in \mathcal{B} \\ l'\in M_{i-1} \text{ ending at $B''$}  \\ A\in \mathcal{A} \, : A\succ B'' \, ; \, A\cap l'=B'' \\ B'\in\mathcal{B} \, ; \, B'\prec A \text{ and } B'\not\in l'}} -\ln\left(\frac{\mu_{i-1}(l')}{\deg(B'').k}\right)\mu_i(l)\\
&=& \mathbb{P}_i\ln(k)+ \sum_{\substack{B''\in \mathcal{B} \\ l'\in M_{i-1} \text{ ending at $B''$}  \\ A\in \mathcal{A} \, : A\succ B'' \, ; \, A\cap l'=B'' \\ B'\in\mathcal{B} \, ; \, B'\prec A \text{ and } B'\not\in l'}} -\ln\left(\frac{\mu_{i-1}(l')}{\deg(B'')}\right)\mu_i(l)\\
&=& \mathbb{P}_i\ln(k)+ \sum_{\substack{B''\in \mathcal{B} \\ l'\in M_{i-1} \text{ ending at $B''$}  \\ A\in \mathcal{A} \, : A\succ B'' \, ; \, A\cap l'=B'' \\ B'\in\mathcal{B} \, ; \, B'\prec A \text{ and } B'\not\in l'}} -\ln\left(\frac{\mu_{i-1}(l')}{\deg(B'')}\right)\frac{\mu_{i-1}(l')}{\deg(B'')k}\\ 
&\geq & \mathbb{P}_i\ln(k)+ \frac{k-i}{k}\sum_{\substack{B''\in \mathcal{B} \\ l'\in M_{i-1} \text{ ending at $B''$}  \\ A\in \mathcal{A} \, : A\succ B'' \, ; \, A\cap l'=B'' }} -\ln\left(\frac{\mu_{i-1}(l')}{\deg(B'')}\right)\frac{\mu_{i-1}(l')}{\deg(B'')}\\
&\geq & \mathbb{P}_i\ln(k) + \frac{k-i}{k}\sum_{\substack{B''\in \mathcal{B} \\ l'\in M_{i-1} \text{ ending at $B''$}}} -\ln\left(\frac{\mu_{i-1}(l')}{\deg(B'')}\right)\mu_{i-1}(l')\frac{\deg(B'')-i+1}{\deg(B'')}
\end{eqnarray*}
\begin{eqnarray*}
&=& \mathbb{P}_i\ln(k) + \frac{k-i}{k}\left[\sum_{\substack{B''\in \mathcal{B} \\ l'\in M_{i-1} \\ \text{ ending at $B''$}}} -\ln\left(\frac{\mu_{i-1}(l')}{\deg(B'')}\right)\mu_{i-1}(l')\right] - \frac{(k-i)(i-1)}{k}\left[\sum_{\substack{B''\in \mathcal{B} \\ l'\in M_{i-1} \\ \text{ ending at $B''$}}} -\ln\left(\frac{\mu_{i-1}(l')}{\deg(B'')}\right)\frac{\mu_{i-1}(l')}{\deg(B'')}\right] \\
&=& \mathbb{P}_i\ln(k) + \frac{k-i}{k}\left[\sum_{\substack{B''\in \mathcal{B} \\ l'\in M_{i-1} \\ \text{ ending at $B''$}}} -\ln(\mu_{i-1}(l'))\mu_{i-1}(l')\right] +  \frac{k-i}{k}\left[\sum_{\substack{B''\in \mathcal{B} \\ l'\in M_{i-1} \\ \text{ ending at $B''$}}}\ln\left(\frac{\deg(B'')}{ka}\right)\mu_{i-1}(l')\right]\\
& & +\frac{k-i}{k}\left[\sum_{\substack{B''\in \mathcal{B} \\ l'\in M_{i-1} \\ \text{ ending at $B''$}}}\ln(ka)\mu_{i-1}(l')\right]- \frac{(k-i)(i-1)}{k}\left[\sum_{\substack{B''\in \mathcal{B} \\ l'\in M_{i-1} \\ \text{ ending at $B''$}}} -\ln\left(\frac{\mu_{i-1}(l')}{\deg(B'')}\right)\frac{\mu_{i-1}(l')}{\deg(B'')}\right].
\end{eqnarray*}

There are 5 terms in the previous two lines, which we will simplify separately. 

\paragraph{The first term\\}The first term, $ \mathbb{P}_i\ln(k)$,  is already simple. We leave it as is.

\paragraph{The second term\\} We just have to use the induction hypothesis:

\[\sum_{l'\in M_{i-1}} -\ln(\mu_{i-1}(l'))\mu_{i-1}(l') \geq \frac{(k-1)!}{(k-i)!k^{i-1})}\ln\left(\left(\frac{ak^2}{b}\right)^{i-1}b\right) - O\left(\frac{b\ln(b)}{ka}\right). \]

\paragraph{The third term\\}
Ignoring the factor of $\frac{k-i}{k}$ at the front for the moment, the third term is $\sum_{\substack{B''\in \mathcal{B} \\ l'\in M_{i-1} \text{ ending at $B''$}}}\ln\left(\frac{\deg(B'')}{ka}\right)\mu_{i-1}(l')$. To simplify this, we will first prove by induction that for all $B\in \mathcal{B}$, we have $\sum_{\substack{l\in M_{i} \\ \text{ starting at $B$}}}\mu_{i}(l) \leq \frac{(k-1)!}{(k-i-1)!k^{i}}\frac{\deg(B)}{ka}$. For $i=1$, we get $\deg(B)(k-1)\frac{1}{ak^2} \leq \frac{(k-1)!}{(k-2)!k^1}\frac{\deg(B)}{ka}$ so it is true for $i=1$. For larger $i$, we will have to remember that we assumed without loss of generality that when picking a new end-vertex $B'$ at step $i$, there are always exactly $k-i$ candidates. From this we get:

\begin{eqnarray*}
& & \sum_{\substack{l\in M_{i} \text{ starting at $B$}}}\mu_{i}(l) \\
& = & \sum_{\substack{l'\in M_{i-1}  \text{ starting at $B$}  \text{ and ending at $B''$} \\ A\in \mathcal{A} \, : A\succ B'' \, ; \, A\cap l'=B'' \\ B'\in\mathcal{B} \, ; \, B'\prec A \text{ and } B'\not\in l' }}\mu_{i-1}(l')\frac{1}{\deg(B')k} \\
& = & \sum_{\substack{l'\in M_{i-1}  \text{ starting at $B$}  \text{ and ending at $B''$} \\ A\in \mathcal{A} \, : A\succ B'' \, ; \, A\cap l'=B'' }}\mu_{i-1}(l')\frac{(k-i)}{\deg(B'')k} \\
&\leq & \sum_{\substack{l'\in M_{i-1}  \text{ starting at $B$}  \text{ and ending at $B''$}  }}\mu_{i-1}(l')\frac{(k-i)\deg(B'')}{\deg(B'')k}\\
& \leq & \left[\frac{(k-1)!}{(k-i)!k^{i-1}}\frac{\deg(B)}{ka}\right]\frac{(k-i)}{k} \\
&=& \frac{(k-1)!}{(k-i-1)!k^{i}}\frac{\deg(B)}{ka} .
\end{eqnarray*}

So by induction, we have proved that for all $B\in \mathcal{B}$, we have $\sum_{\substack{l\in M_{i} \\ \text{ starting at $B$}}}\mu_{i}(l) \leq \frac{(k-1)!}{(k-i-1)!k^{i}}\frac{\deg(B)}{ka}$. Now notice that $\mu_i$ is symmetrical with respect to start and finish (remember that if we write the vertices of the path as $B_0$, $B_1$, ..., $B_i$, then $\mu_i(l)=\frac{1}{\deg(B_1)\deg(B_2)...\deg(B_{i-1})ak^{i+1}}$). Therefore we get that for all $B\in \mathcal{B}$, $\sum_{\substack{l\in M_{i} \\ \text{ ending at $B$}}}\mu_{i}(l) \leq \frac{(k-1)!}{(k-i-1)!k^{i}}\frac{\deg(B)}{ka}$. Now using this result, we can go back to trying to estimate the third term:

\begin{eqnarray*}
& &\sum_{\substack{B''\in \mathcal{B} \\ l'\in M_{i-1} \text{ ending at $B''$}}}\ln\left(\frac{\deg(B'')}{ka}\right)\mu_{i-1}(l')\\
&\geq & \sum_{B''\in \mathcal{B}}\ln\left(\frac{\deg(B'')}{ka}\right)\frac{(k-1)!}{(k-i)!k^{i-1}}\frac{\deg(B'')}{ka}\\
&=& \frac{(k-1)!}{(k-i)!k^{i-1}}\sum_{B''\in \mathcal{B}}\ln\left(\frac{\deg(B'')}{ka}\right)\frac{\deg(B'')}{ka}\\
&=& \frac{(k-1)!}{(k-i)!k^{i-1}}\sum_{B''\in \mathcal{B}}\ln\left(\mu_0(B'')\right)\mu_0(B'').
\end{eqnarray*}

Using the entropy inequality on $\mu_0$ again, we get that this is at least:

\[ \frac{(k-1)!}{(k-i)!k^{i-1}}\sum_{B''\in \mathcal{B}}\ln\left(\frac{1}{b}\right)\frac{1}{b} = -\frac{(k-1)!}{(k-i)!k^{i-1}}\ln(b).\]

\paragraph{The fourth term\\} We have $\sum_{l'\in M_{i-1}}\ln(ka)\mu_{i-1}(l')=\ln(ka)\mathbb{P}_{i-1}$ by definition of $\mathbb{P}_{i-1}$. 

\paragraph{The fifth term\\}

\[\sum_{\substack{B''\in \mathcal{B} \\ l'\in M_{i-1} \text{ ending at $B''$}}}\ln\left(\frac{\mu_{i-1}(l')}{\deg(B'')}\right)\frac{\mu_{i-1}(l')}{\deg(B'')} \geq \sum_{\substack{B''\in \mathcal{B} \\ l'\in L_{i-1} \text{ ending at $B''$}}}\ln\left(\frac{\mu_{i-1}(l')}{\deg(B'')}\right)\frac{\mu_{i-1}(l')}{\deg(B'')}.\]

Using property (\ref{prop1}) again, we get that $\sum_{\substack{B''\in\mathcal{B} \\ l'\in L_{i-1} \text{ ending at $B''$}}}\frac{\mu_{i-1}(l')}{\deg(B'')} = \sum_{B''\in \mathcal{B}}\frac{1}{ka}=\frac{b}{ka}$. We can now do another entropy inequality to get that the fifth term is at least: 

\begin{eqnarray*}
& &\sum_{\substack{B''\in \mathcal{B} \\ l'\in L_{i-1} \text{ ending at $B''$}}}\ln\left(\frac{\mu_{i-1}(l')}{\deg(B'')}\right)\frac{\mu_{i-1}(l')}{\deg(B'')}\\
& = & \frac{b}{ka}\left(\ln(b/ka) + \sum_{\substack{B''\in \mathcal{B} \\ l'\in L_{i-1} \text{ ending at $B''$}}}\ln\left(\frac{\mu_{i-1}(l')\cdot ka}{\deg(B'')b}\right)\frac{\mu_{i-1}(l')\cdot ka}{\deg(B'')b}\right) \\
& \geq &  \frac{b}{ka}\left(\ln(b/ka) + \sum_{\substack{B''\in \mathcal{B} \\ l'\in L_{i-1} \text{ ending at $B''$}}}\ln\left(\frac{1}{|L_{i-1}|}\right)\frac{1}{|L_{i-1}|}\right) \\
&= & \sum_{\substack{B''\in \mathcal{B} \\ l'\in L_{i-1} \text{ ending at $B''$}}}\ln\left(\frac{b}{ka|L_{i-1}|}\right)\frac{b}{ka|L_{i-1}|}\\
&=& \ln\left(\frac{b}{ka|L_{i-1}|}\right)\frac{b}{ka} .
\end{eqnarray*}

We now need to bound $L_{i-1}$ from above. We don't need to do anything fancy for this. We start by picking an edge with 2 vertices on it (there are $ak^2$ possibilities). Then pick an edge incident to the ending vertex. There are at most $b$ possible choices. Choose a new ending vertex somewhere on the new edge; there are $k$ possibilities. Repeat this process $i-2$ times until we have an element of $L_{i-1}$. This gives us $|L_{i-1}| \leq ak^ib^{i-2}$. Therefore the fifth term is at least:

\[ \geq \ln\left(\frac{1}{a^2k^{i+1}b^{i-2}}\right)\frac{b}{ka} = O\left(\frac{b\ln(b)}{ka}\right).\]

\paragraph{Putting the five terms back together again\\}

Thus, by adding up the five terms back together again, we get:

\begin{eqnarray*} 
& & \sum_{l\in M_i} -\ln(\mu_i(l))\mu_i(l)\\
& \geq & \mathbb{P}_i\ln(k)+\frac{k-i}{k}\left[\frac{(k-1)!}{(k-i)!k^{i-1})}\ln\left(\left(\frac{ak^2}{b}\right)^{i-1}b\right) - O\left(\frac{b\ln(b)}{ka}\right)\right]\\
& & +\frac{k-i}{k}\left[-\frac{(k-1)!}{(k-i)!k^{i-1}}\ln(b)\right]+\frac{k-i}{k}\left[\ln(ka)\mathbb{P}_{i-1}\right]-O\left(\frac{b\ln(b)}{ka}\right).
\end{eqnarray*}

Now using the fact from (\ref{prob}) that $\mathbb{P}_{i-1}\geq \frac{(k-1)!}{(k-i)!k^{i-1}} -O\left(\frac{b}{ka}\right) $ and $\mathbb{P}_{i}\geq \frac{(k-1)!}{(k-i-1)!k^{i}} -O\left(\frac{b}{ka}\right) $, we get that it is at least

\begin{eqnarray*}
& & \frac{(k-1)!}{(k-i-1)!k^{i}} \left[\ln(k)+\ln\left(\left(\frac{ak^2}{b}\right)^{i-1}b\right)-\ln(b)+\ln(ka)\right] -O\left(\frac{b\ln(b)}{ka}\right) \\
&=&\frac{(k-1)!}{(k-i-1)!k^{i}} \ln\left(\left(\frac{ak^2}{b}\right)^{i}b\right) -O\left(\frac{b\ln(b)}{ka}\right).
\end{eqnarray*}

This completes the proof of proposition (\ref{prop2}) by induction.

\paragraph{Finishing up the lower bound on $|M_i|$\\}

Now we can use (\ref{prop3}), which tells us that $\ln(|M_i|) \geq \ln(\mathbb{P}_i)+\frac{(k-1)!}{(k-i-1)!k^i\mathbb{P}_i}\ln\left(\frac{ak^2}{b}\right) - O\left(\frac{b\ln(b)}{ka}\right)$. Using what we know about the probability $\mathbb{P}_i$ from (\ref{prob}), we get the following:

\begin{eqnarray*}
 \ln(|M_i|) & \geq & \ln\left(\frac{(k-1)!}{(k-i-1)!k^{i}}\left(1-O\left(\frac{b}{ka}\right)\right)\right) +\frac{\frac{(k-1)!}{(k-i-1)!k^{i}} \ln\left(\left(\frac{ak^2}{b}\right)^{i}b\right) -O\left(\frac{b\ln(b)}{ka}\right)}{\frac{(k-1)!}{(k-i-1)!k^{i}}\left(1-O\left(\frac{b}{ka}\right)\right)}\\
&=& \ln\left(\frac{(k-1)!}{(k-i-1)!k^{i}}\left(1-O\left(\frac{b}{ka}\right)\right)\right) + \ln\left(\left(\frac{ak^2}{b}\right)^{i}b\right) - O\left(\frac{b\ln(b)}{ka}\right)\\
&=& \ln\left(\frac{(k-1)!}{(k-i-1)!} \left(\frac{ak}{b}\right)^{i}b\left(1-O\left(\frac{b\ln(b)}{ka}\right)\right)\right).
\end{eqnarray*}

Thus: 

\begin{equation}
|M_i|\geq \frac{(k-1)!}{(k-i-1)!} \left(\frac{ak}{b}\right)^{i}b\left(1-O\left(\frac{b\ln(b)}{ka}\right)\right). \label{M_i}\end{equation}

\subsection{Upper bound on $a$ as a function of $b$}

Remember that $M_i$ was the set of paths between pairs at distance $i$ from each other (with direction). Given any such pair, say $B$ and $B'=B\cup\{x_1,x_2,...,x_{i}\}\backslash\{y_1,y_2,...,y_{i}\}$ there are at most $i!^2$ such paths. This is because to get a path of length $i$, you need to pick some sequence of the $x_j$s and $y_j$s and you then follow the path by adding the next $x_j$ and subtracting the next $y_j$ at each step; this then determines the path uniquely. So therefore the total number of pairs, (which is equal to $b(b-1)$ ), is at least $\sum_{i=1}^{k-1} \frac{|M_i|}{i!^2}$. Thus:

\[b(b-1)\geq \sum_{i=1}^{k-1}\left(\frac{ak}{b}\right)^i\frac{(k-1)!}{(k-i-1)!i!^2}b\left(1-O\left(\frac{b\ln(b)}{ka}\right)\right)\]

\[b-1 \geq \left(\frac{ak}{b}\right)^{k-1}/(k-1)! + O\left(\left(\frac{ak}{b}\ln(b)\right)^{k-2}\right)\]

\[\frac{ak}{b} \leq (b(k-1)!)^{1/(k-1)} + O(\ln(b))\]

\[a \leq b^{k/(k-1)}\frac{(k-1)!^{1/(k-1)}}{k} + O(b\ln(b).)\]

This matches our canonical example where $a=\binom{c}{k}$ and $b=\binom{c}{k-1}$ to within $O(b\ln(b))$. Indeed, the canonical example has $\frac{ak}{b}=c-k+1$ and $b(k-1)!=c(c-1)...(c-k+2) = c^{k-1}+O(c^{k-2})$, so $(b(k-1)!)^{1/(k-1)}=c+O(1)$, therefore $\frac{ak}{b} = (b(k-1)!)^{1/(k-1)} +O(1)$ in our canonical example. This matches our result to within $O(b\ln(b))$.

\subsection{Using stability to gather information about our sets}

This is similar to what we did in the case $k=3$. Suppose we have a valid configuration with $b$ vertices and $a = b^{k/(k-1)}\frac{(k-1)!^{1/(k-1)}}{k}\left(1+O\left(\frac{\ln(b)}{b^{1/{k-1}}}\right)\right)$ edges. We will go through the proof to see what properties we can deduce of $\mathcal{A}$ and $\mathcal{B}$.

\paragraph{Nice pairs\\} This is based on Section 6.5 of the original proof. We know from (\ref{M_i}) that \[|M_{k-1}|\geq (k-1)! \left(\frac{ak}{b}\right)^{k-1}b\left(1-O\left(\frac{b\ln(b)}{ka}\right)\right) = (k-1)!^2 b^2\left(1+O\left(\frac{\ln(b)}{b^{1/(k-1)}}\right)\right).\]
We know that there are less than $b^2$ pairs of points at distance $k-1$. The maximum number of paths in $M_{k-1}$ joining such a pair is $(k-1)!^2$. Most pairs should have exactly $(k-1)!^2$ paths; we'll call this a \underline{nice pair}. However, there might be some that have less. We'll say that there are $\omega$ nice pairs. Then the number of paths of length $k-1$ is less than $b^2(k-1)!^2 - \omega$. Combining this with our inequality for $M_{k-1}$, this tells us that $\omega\leq O\left(b^2\frac{\ln(b)}{b^{1/(k-1)}}\right)$. So almost all pairs at distance $k-1$ will be nice. The exceptions only make up a proportion of $O\left(\frac{\ln(b)}{b^{1/(k-1)}}\right)$ of the total. We also get that $|M_{k-1}|=(k-1)!^2 b^2\left(1+O\left(\frac{\ln(b)}{b^{1/(k-1)}}\right)\right)$.\\
\\
\\
So given a vertex $B$, the average number of points $B''$ that create a nice pair with it is $b\left(1-O\left(\frac{\ln(b)}{b^{1/(k-1)}}\right)\right)$. We want to know how many vertices are close to that number. To be more precise, lets say that there are $b(1-\delta)$ vertices $B$ that have at least $b(1-\epsilon)$ vertices $B''$ at distance $k-1$ from it (where $\epsilon$ and $\delta$ will be defined later). Then the total number of pairs at distance $k-1$ is:

\[b^2\left(1-O\left(\frac{\ln(b)}{b^{1/(k-1)}}\right)\right) \leq b(b-b\epsilon) + (b-b\delta)(b\epsilon).\]

Reordering the inequality gives:

\[\epsilon \leq \frac{1}{\delta}O\left( \frac{\ln(b)}{b^{1/(k-1)}}\right).\]

So if we set $\delta=b^{-1/(2k-2)}$ and $\epsilon=O\left( \frac{\ln(b)}{b^{1/(2k-2)}}\right)$, then this works. Therefore there are $b\left(1-b^{-1/(2k-2)}\right)$ vertices $B$ that are part of at least $b\left(1-O\left( \frac{\ln(b)}{b^{1/(2k-2)}}\right)\right)$ nice pairs. We'll call this set of vertices $\mathcal{B'}$.
 
\paragraph{Low degree\\}

We know that the sum of the degrees of all the vertices is $ka$. Therefore the sum of the degrees of vertices in $\mathcal{B'}$ is at most $ka$. So the average degree of an element of $\mathcal{B'}$ is at most $\frac{ka}{b\left(1-b^{-1/(2k-2)}\right)} \leq (b(k-1)!)^{1/(k-1)}(1+\epsilon)$. \\
\\
Therefore there exists a vertex $B$ of $\mathcal{B'}$ with degree less than $(b(k-1)!)^{1/(k-1)}(1+\epsilon).$

\paragraph{Colouring the edges incident to $B$\\}

We know we have a vertex $B$ that is part of at least $b(1-\epsilon)$ nice pairs, and moreover, it has degree $d\leq (b(k-1)!)^{1/(k-1)}(1+\epsilon)$. We'll denote the set of edges incident to $B$ by $\mathcal{G} = \{B\cup\{x_1\}, B\cup\{x_2\}, ..., B\cup\{x_d\}\}$. \\
\\
Each of the edges in $\mathcal{G}$ has $k-1$ other vertices incident to it: $B\cup\{x_i\}\backslash\{y_{i,1}\}$, $B\cup\{x_i\}\backslash\{y_{i,2}\}$, ..., $B\cup\{x_i\}\backslash\{y_{i,k-1}\}$. We'll also colour $\mathcal{G}$ by giving colour: $\{y_{i,1},y_{i,2},...,y_{i,k-1}\}$ to the edge $B\cup\{x_i\}$.\\
\\
Now for every nice pair $(B,B'')$, there are $(k-1)!^2$ paths between them. $B''$ is at distance $k-1$ from $B$ so write $B''=B\cup\{t_1,t_2,...,t_{k-1}\}\backslash\{w_1,w_2,...,w_{k-1}\}$. Now there exist all $(k-1)!^2$ possible paths between $B$ and $B''$, which means that $B\cup\{t_1\}$, $B\cup\{t_2\}$, ..., $B\cup\{t_{k-1}\}$ are all in $\mathcal{G}$. Therefore , each $t_i$ is equal to to some $x_j$. Furthermore, all these edges have the same colour: $\{w_1,w_2,...,w_{k-1}\}$. So therefore we know that for every nice pair $(B,B'')$, there exists a corresponding set of $k-1$ elements of $\mathcal{G}$ that all have the same colour. Furthermore, given such a monochromatic set of $k-1$ elements of $\mathcal{G}$, there is at most one $B''$ that they correspond to. So the number of nice pairs containing $B$ is at most the number of monochromatic $(k-1)$-sets of $\mathcal{G}$.\\
\\
We'll say that the largest colour class of $\mathcal{G}$ has size $d(1-\alpha)$. Then the maximum number of monochromatic $(k-1)$-sets is $\binom{d\alpha}{k-1}+\binom{d(1-\alpha)}{k-1} = \frac{d^{k-1}}{(k-1)!}(\alpha^{k-1}+(1-\alpha)^{k-1}-O(1/d))$. We know that this number is at least $b(1-\epsilon)$ and so we plug in $d\leq (b(k-1)!)^{1/(k-1)}(1+\epsilon)$ to get:

\begin{eqnarray*}
\frac{b(k-1)!}{(k-1)!}(\alpha^{k-1}+(1-\alpha)^{k-1}-O(1/d))(1+\epsilon)^{k-1} &\geq & b(1-\epsilon)\\
(\alpha^{k-1}+(1-\alpha)^{k-1}-O(1/d)) & \geq & \frac{1-\epsilon}{(1+\epsilon)^{k-1}} \\
\alpha & \leq & \frac{k}{k-1}\epsilon + O(\epsilon^2)+O(1/d).
\end{eqnarray*}

Therefore, we know that there is a very large colour class of size $d(1-\frac{k\epsilon}{k-1}+ O(\epsilon^2)+O(1/d))$, comprising nearly all elements of $\mathcal{G}$. We'll say its colour is $\{z_1,z_2,...,z_{k-1}\}$.

\paragraph{A nice hypergraph\\}

We want to know how many elements of $\mathcal{B'}$ are connected to our large colour class. The maximum number of them that we don't use in our large colour class is $\binom{\alpha}{k-1} \leq \left(\frac{k\epsilon}{k-1}\right)^{k-1}\frac{d^{k-1}}{(k-1)!} \leq \left(\frac{k\epsilon(1+\epsilon)}{k-1}\right)^{k-1}b$. Therefore the number of vertices of $\mathcal{B'}$ that are connected to our large colour class is at least $b\left(1-\epsilon-\left(\frac{k\epsilon(1+\epsilon)}{k-1}\right)^{k-1}\right)$, so that is nearly all points.\\
\\
Now notice that every one of these vertices of $\mathcal{B'}$ that is connected to our large colour class contains $B\backslash\{z_1,z_2,...,z_{k-1}\}$ because that is the only way to connect it to edges in $\mathcal{G}$ of that colour. We set $S=B\backslash\{z_1,z_2,...,z_{k-1}\}$ and we end up with a nice hypergraph that comprises nearly all the vertices of $\mathcal{B}$. So what we have is:

\begin{mylem}
There exists a set $S\in \mathbb{N}^{(r-k)}$ such that $S$ is a subset of $b(1-o(1))$ elements of $\mathcal{B}$.
\end{mylem}

We'll define our nice hypergraph $\mathcal{D}$ to consist of all the vertices an edges that contain $S$ as a subset.

\subsection{Using classical Kruskal Katona to improve the bound further}

At this point, we know that we have a large nice hypergraph of vertices $\mathcal{D}$, all of which contain $S$ as a subset. We'll say that there are $\lambda$ vertices in $\mathcal{B}\backslash\mathcal{D}$. We know that $\lambda=o(b)$. Our aim in this section is to bound $\lambda$ by a constant. How many edges can we have in our graph? To count them, we'll separate them into 3 cases:\\
\\
$\bullet$ The edges that are entirely contained within $\mathcal{D}$. We can apply the classical version of the Kruskal Katona Theorem to get an upper bound. To get that bound, we need to write $|\mathcal{D}| = b-\lambda$ in the form $\binom{d_{k-1}}{k-1} + \binom{d_{k-2}}{k-2} + ... + \binom{d_{1}}{1}$. Then the maximum number of edges is $\binom{d_{k-1}}{k} + \binom{d_{k-2}}{k-1} + ... + \binom{d_{1}}{2}$.\\
\\
$\bullet$ The edges that are entirely contained within $\mathcal{B}\backslash\mathcal{D}$. For these, we can just apply our formula to say that there are at most $\lambda^{k/(k-1)}\frac{(k-1)!^{1/(k-1)}}{k}\left(1+O\left(\frac{\ln(\lambda)}{\lambda^{1/(k-1)}}\right)\right)$ of them.\\
\\
$\bullet$ The edges that are incident to both $\mathcal{D}$ and $\mathcal{B}\backslash\mathcal{D}$. Since these edges are incident to some vertex in our nice hypergraph $\mathcal{D}$, that vertex has to contain $S$ as a subset, therefore the edge also has to contain $S$. Now for every vertex $B$ in $\mathcal{B}\backslash\mathcal{D}$, $S\not\subset B$, so the only edge that can connect $B$ to $\mathcal{D}$ has to be $B\cup S$. In particular, it is unique. Therefore the number of edges that fall under this case is at most $\lambda$. \\
\\
If we add up everything, we get that the maximal number of edges is:\\ $\binom{d_{k-1}}{k} + \binom{d_{k-2}}{k-1} + ... + \binom{d_{1}}{2} + O(\lambda^{k/(k-1)})$, where $\lambda=b-\left[\binom{d_{k-1}}{k-1} + \binom{d_{k-2}}{k-2} + ... + \binom{d_{1}}{1}\right]$. \\
\\
We want to compare this to what we would get with our hypothesis (which states that $\lambda=0$ is optimal). For this, you would write $b=\binom{c_{k-1}}{k-1} + \binom{c_{k-2}}{k-2} + ... + \binom{c_{1}}{1}$, and then the number of edges would be: $\binom{c_{k-1}}{k} + \binom{c_{k-2}}{k-1} + ... + \binom{c_{1}}{2}$. So if we did have a counter-example to our hypothesis, we would have:

\[\left[\binom{c_{k-1}}{k} + \binom{c_{k-2}}{k-1} + ... + \binom{c_{1}}{2}\right] - \left[\binom{d_{k-1}}{k} + \binom{d_{k-2}}{k-1} + ... + \binom{d_{1}}{2}\right] = O(\lambda^{k/(k-1)})\]

\[\text{where } \lambda = \left[\binom{c_{k-1}}{k-1} + \binom{c_{k-2}}{k-2} + ... + \binom{c_{1}}{1}\right] - \left[\binom{d_{k-1}}{k-1} + \binom{d_{k-2}}{k-2} + ... + \binom{d_{1}}{1}\right]\].

\paragraph{Case 1: $c_{k-1} -d_{k-1} \geq 2$\\}

Then $\lambda \leq \binom{c_{k-1}+1}{k-1}-\binom{d_{k-1}}{k-1} \leq (c_{k-1}+1-d_{k-1})\binom{c_{k-1}}{k-2} = (c_{k-1}+1-d_{k-1}) O((c_{k-1})^{k-2})$. Therefore the right hand side of the inequality is at most $(c_{k-1}+1-d_{k-1})^{k/(k-1)}O\left((c_{k-1})^{k(k-2)/(k-1)}\right)$. \\
\\
Meanwhile, the left hand side of the inequality is at least $\binom{c_{k-1}}{k}-\binom{d_{k-1}+1}{k} \geq (c_{k-1}-d_{k-1}-1)\binom{d_{k-1}+1}{k-1} = (c_{k-1}-d_{k-1}-1)\Omega(c_{k-1}^{k-1})$. \\
\\
Notice that $(c_{k-1}-d_{k-1}-1) \geq 1$, so to get the inequality to hold, we must have $\frac{(c_{k-1}-d_{k-1}+1)^{k/(k-1)}}{c_{k-1}-d_{k-1}-1} > \Omega((c_{k-1})^{1/(k-1)})$ , so therefore $c_{k-1}-d_{k-1} = \Omega(c_{k-1})$. But we know that $\lambda = o(b)$ so $(c_{k-1})^{k-1}-(d_{k-1})^{k-1} = o((c_{k-1})^{k-1})$ so $c_{k-1}-d_{k-1} = o(c_{k-1})$. This is a contradiction, so the inequality never holds when $c_{k-1} -d_{k-1} \geq 2$, so there are no counter-examples of this type (as long as $c_{k-1}$ is large).

\paragraph{Case 2: $c_{k-1} -d_{k-1} =1$\\}

We substitute $d_{k-1} = c_{k-1}-1$ into the inequality to get:

\[\binom{c_{k-1}-1}{k-1} + \left[\binom{c_{k-2}}{k-1} + ... + \binom{c_{1}}{2}\right] - \left[\binom{d_{k-2}}{k-1} + ... + \binom{d_{1}}{2}\right] < O(\lambda^{k/(k-1)})\]

\[\text{where } \lambda = \binom{c_{k-1}-1}{k-2} + \left[ \binom{c_{k-2}}{k-2} + ... + \binom{c_{1}}{1}\right] - \left[\binom{d_{k-2}}{k-2} + ... + \binom{d_{1}}{1}\right].\]

Then $\lambda \leq \binom{c_{k-1}}{k-2} = O((c_{k-1})^{k-2})$. Therefore the right hand side is at most $O((c_{k-1})^{k(k-2)/(k-1)})$.\\
\\
Meanwhile, the left hand side of the inequality is at least:\\ $  \binom{c_{k-1}-1}{k-1} - \binom{d_{k-2}+1}{k-1}  = \\ \frac{(c_{k-1})^{k-1}}{(k-1)!} - O((c_{k-1})^{k-2})- \frac{(d_{k-2})^{k-1}}{(k-1)!} + O((d_{k-2})^{k-2})$. The only way to get this to be smaller than the right hand side is to have the $\frac{(c_{k-1})^{k-1}}{(k-1)!} - \frac{(d_{k-2})^{k-1}}{(k-1)!} = O((c_{k-1})^{k-1-1/(k-1)})$ so $c_{k-1}-d_{k-2} = O((c_{k-1})^{1-1/(k-1)})) = o(c_{k-1})$.\\
\\
Now we try again except this time we can use the information that $c_{k-1}-d_{k-2} = o(c_{k-1})$. We have $\lambda \leq \binom{c_{k-1}-1}{k-2} + \binom{c_{k-2}+1}{k-2} - \binom{d_{k-2}}{k-2}  \leq (c_{k-1}-d_{k-2}-1)\binom{c_{k-1}-2}{k-3} + \binom{c_{k-2}+1}{k-2}$. We can use Jensen's inequality to deduce that the right hand side of the inequality is at most:\\ $(c_{k-1}-d_{k-2}-1)^{k/(k-1)}O((c_{k-1})^{k(k-3)/(k-1)}) + O((c_{k-2})^{k(k-2)/(k-1)})$. \\
\\
Meanwhile, the left hand side of the inequality is at least: $\binom{c_{k-1}-1}{k-1} + \binom{c_{k-2}}{k-1} - \binom{d_{k-2}+1}{k-1} \geq \\ (c_{k-1}-d_{k-2}-2)\binom{d_{k-2}+1}{k-2} + \binom{c_{k-2}}{k-1} = (c_{k-1}-d_{k-2}-2)\Omega((c_{k-1})^{k-2}) + \Omega((c_{k-2})^{k-1})$. The only way this is smaller than the right hand side is if $(c_{k-1}-d_{k-2}-2)\Omega((c_{k-1})^{k-2}) < (c_{k-1}-d_{k-2}-1)^{k/(k-1)}O((c_{k-1})^{k(k-3)/(k-1)})$. This implies that $\frac{c_{k-1}-d_{k-2}-2}{(c_{k-1}-d_{k-2}-1)^{k/(k-1)}} < O((c_{k-1})^{-2/(k-1)})$.\\
\\
There are two solutions to this: either $c_{k-1}-d_{k-2}-2=0$, or $c_{k-1}-d_{k-2} = \Omega((c_{k-1})^2)$. But the second solution is clearly impossible, so the only possibility is $d_{k-2}=c_{k-1}-2$. We substitute this back into the inequality and we get:

\[\binom{c_{k-1}-2}{k-2} + \left[\binom{c_{k-2}}{k-1} + ... + \binom{c_{1}}{2}\right] - \left[\binom{d_{k-3}}{k-2} + ... + \binom{d_{1}}{2}\right] < O(\lambda^{k/(k-1)})\]

\[\text{where } \lambda = \binom{c_{k-1}-2}{k-3} + \left[ \binom{c_{k-2}}{k-2} + ... + \binom{c_{1}}{1}\right] - \left[\binom{d_{k-3}}{k-3} + ... + \binom{d_{1}}{1}\right].\]

This looks remarkably like the inequality we had at the start of this case. In fact, we can repeat the argument with a few changes to prove that $d_{k-3}=c_{k-1}-3$. Then we can continue using the same argument to prove $d_{k-4}=c_{k-1}-4$, ... all the way until $d_{1}=c_{k-1}-(k-1)$. At this point, we're left with:

\[(c_{k-1}-k+1) + \left[\binom{c_{k-2}}{k-1} + ... + \binom{c_{1}}{2}\right]  < O(\lambda^{k/(k-1)})\]

\[\text{where } \lambda = 1 + \left[ \binom{c_{k-2}}{k-2} + ... + \binom{c_{1}}{1}\right].\]

The we need $\lambda^{k/(k-1)} \geq \Omega(c_{k-1})$ so $\binom{c_{k-2}+1}{k-2} \geq \lambda \geq \Omega((c_{k-1})^{(k-1)/k})$ so $c_{k-2} \geq \Omega((c_{k-1})^{(k-1)/k/(k-2)})$. Because $\frac{(k-1)^2}{k(k-2)} >1$, the dominant term on the left hand side of the inequality is $\binom{c_{k-2}}{k-1}$. Therefore we end up with $\Omega((c_{k-2})^{k-1}) < O((c_{k-2})^{k(k-2)/(k-1)})$ which is impossible for $c_{k-2}$ large enough. $c_{k-2} \geq \Omega((c_{k-1})^{(k-1)/k/(k-2)})$ so it's also impossible when $c_{k-1}$ is large enough. So the inequality never holds and there are no large counter-examples of this type.

\paragraph{Case 3: $c_{k-1}-d_{k-1}=0$\\}

We substitute $d_{k-1}=c_{k-1}$ into the inequality to get:

\[\left[\binom{c_{k-2}}{k-1} + ... + \binom{c_{1}}{2}\right] - \left[\binom{d_{k-2}}{k-1} + ... + \binom{d_{1}}{2}\right] < O(\lambda^{k/(k-1)})\]

\[\text{where } \lambda = \left[ \binom{c_{k-2}}{k-2} + ... + \binom{c_{1}}{1}\right] - \left[\binom{d_{k-2}}{k-2} + ... + \binom{d_{1}}{1}\right].\]

But $O(\lambda^{k/(k-1)}) < O(\lambda^{(k-1)/(k-2)})$, so we get: 

\[\left[\binom{c_{k-2}}{k-1} + ... + \binom{c_{1}}{2}\right] - \left[\binom{d_{k-2}}{k-1} + ... + \binom{d_{1}}{2}\right] < O(\lambda^{(k-1)/(k-2)})\]

\[\text{where } \lambda = \left[ \binom{c_{k-2}}{k-2} + ... + \binom{c_{1}}{1}\right] - \left[\binom{d_{k-2}}{k-2} + ... + \binom{d_{1}}{1}\right].\]

This new inequality is identical to our original except that we have $k-1$ instead of $k$. If $c_{k-2}$ is large enough, then we can repeat our argument and either go to cases 1 and 2 and prove no counter-example exists, or go back through case 3 and reduce $k$ again (more formally, there exists a constant $\mu$ depending only on $k$ such that if for any $i$,  $d_i<c_i>\mu$, then case 1 or 2 applies and there is no counter-example possible). If all of $c_{k-2}$, $c_{k-3}$, ..., $c_{1}$ are larger than $\mu$, then we get $c_{k-1}=d_{k-1}$, $c_{k-2}=d_{k-2}$, ... , $c_{1}=d_{1}$, which implies that $b=b-\lambda$ so $\lambda=0$.\\
\\
\\
This is exactly Theorem 2: there is some constant $\mu$ depending only on $k$ such that if $b=\left[ \binom{c_{k-1}}{k-1} + ... + \binom{c_{1}}{1}\right]$, for some $c_{k-1}>c_{k-2}>...>c_1>\mu$, then the maximum value for $a$ is exactly:

\[\left[ \binom{c_{k-1}}{k} + ... + \binom{c_{1}}{2}\right] = f(r,k,b).\]
The only cases that aren't covered by Theorem 2, is if there exists some $i$ such that $c_{i}$ is smaller than $\mu$. Without loss of generality, pick $i$ to be the largest such. Then we know that $c_{j}=d_{j}$ for any $j>i$. Then we get:

\[\lambda \leq  \left[ \binom{\mu}{i} \binom{\mu-1}{i-1}+ ... + \binom{\mu-i+1}{1}\right] -[0] .\]

So $\lambda$ is in fact bounded by a constant, which implies the following lemma: 

\begin{mylem}
	There is a constant $\lambda_{max}$ depending only on $k$ such that there exists a subset $\mathcal{D}$ of $\mathcal{B}$ of size at least $b-\lambda_{max}$ and a set $S\in \mathbb{N}^{(r-k)}$ such that $S$ is a subset of every element of $\mathcal{D}$.
\end{mylem}

And now Theorem 3 is just an easy corrolary of this, since we have a bounded number of vertices that aren't in our nice hypergraph, these vertices can only form a bounded number of extra edges, therefore there is a constant $\tau$ depending only on $k$ such that if $b=\left[ \binom{c_{k-1}}{k-1} + ... + \binom{c_{1}}{1}\right]$, for some $c_{k-1}>c_{k-2}>...>c_1$, then the maximum value for $a$ is between:
	
	\[\left[ \binom{c_{k-1}}{k} + ... + \binom{c_{1}}{2}\right] \leq f(r,k,b) \leq \left[ \binom{c_{k-1}}{k} + ... + \binom{c_{1}}{2}\right] + \tau.\]

\bibliography{Kruskal_Katona}{}
\bibliographystyle{plain}

\end{document}